\documentclass[11pt,a4paper,mathrsfs,leqno]{article}
\usepackage{amsmath,amssymb,amsxtra,epsf,amscd,enumerate,ulem,graphics,color,fancybox}
\usepackage[english]{babel}
\usepackage[all]{xy}

\newtheorem{theorem}{Theorem}[section]

\newtheorem{proposition}[theorem]{Proposition}

\newtheorem{definition}[theorem]{Definition}

\newtheorem{corollary}[theorem]{Corollary}

\newtheorem{lemma}[theorem]{Lemma}

\newtheorem{claim}[theorem]{Claim}

\newtheorem{conjecture}[theorem]{Conjecture}

\newtheorem{remark}[theorem]{{Remark}}
\newenvironment{rmk}{\begin{remark}\rm }{\end{remark}}
\newtheorem{example}[theorem]{{Example}}
\newenvironment{ex}{\begin{example}\rm }{\end{example}}

\newenvironment{proof}{{\flushleft{\textit{Proof:} }} \rm }{\Qed \\}

\newcommand{\bp}{\begin{proposition}} \newcommand{\Bpo}{\bp}
\newcommand{\ep}{\end{proposition}} \newcommand{\Epo}{\ep}
\newcommand{\bd}{\begin{definition}} \newcommand{\Bdf}{\bd}
\newcommand{\ed}{\end{definition}} \newcommand{\Edf}{\ed} 
\newcommand{\bt}{\begin{theorem}} \newcommand{\Bte}{\bt}
\newcommand{\et}{\end{theorem}} \newcommand{\Ete}{\et}
\newcommand{\bc}{\begin{corollary}} 
\newcommand{\ec}{\end{corollary}} 
\newcommand{\bl}{\begin{lemma}} \newcommand{\Blm}{\bl}
\newcommand{\el}{\end{lemma}} \newcommand{\Elm}{\el}
\newcommand{\br}{\begin{rmk}} 
\newcommand{\er}{\end{rmk}} 
\newcommand{\bex}{\begin{ex}} 
\newcommand{\eex}{\end{ex}}
\newcommand{\bpf}{\begin{proof}} \newcommand{\Bdm}{\bpf}
\newcommand{\epf}{\end{proof}} \newcommand{\Edm}{\epf}

\newcommand{\Bca}{\begin{claim}}
\newcommand{\Eca}{\end{claim}}
\newcommand{\Bcn}{\begin{conjecture}}
\newcommand{\Ecn}{\end{conjecture}}

\numberwithin{equation}{section}

\newcommand{\dsty}{\displaystyle}
\newcommand{\Qed}{\hfill \mbox{$\blacksquare$}}

\newcommand{\ppq}{\leqslant}
\newcommand{\pgq}{\geqslant}
\newcommand{\ot}{\otimes}
\newcommand{\ots}{\otimes_{\kk Q_0}}
\newcommand{\ovl}{\overline}

\newcommand{\ci}{\frak{c}}

\newcommand{\nn}{\mathbb{N}}
\newcommand{\kk}{\Bbbk}

\newcommand{\car}{\mathrm{char}}
\newcommand{\Hom}{\mathrm{Hom}}
\newcommand{\Ext}{\mathrm{Ext}}
\newcommand{\Tor}{\mathrm{Tor}}
\newcommand{\hh}{\mathrm{H}}
\newcommand{\id}{\mathrm{id}}
\newcommand{\im}{\mathrm{Im}}
\renewcommand{\ker}{\mathrm{Ker}}

\newcommand{\C}{\mathcal{C}}
\newcommand{\I}{\mathcal{I}}

\newcommand{\set}[1]{\left\{ #1 \right\}}

\newcommand{\mx}[1]{\begin{pmatrix}#1\end{pmatrix}}

\newcommand{\pot}[1]{\mathrm{Pot}(#1)}
\newcommand{\der}[1]{\frac{\partial}{\partial{#1}}}

\title{\textsc{Poincar\'e-Birkhoff-Witt deformations of Calabi-Yau algebras}}
\author{Roland {\sc  Berger}, Rachel {\sc Taillefer} }
\date{}

\begin{document}

\maketitle

\begin{abstract} Recently, Bocklandt proved a conjecture by Van den Bergh in its graded version,  stating that a graded quiver algebra (with relations) which is Calabi-Yau of dimension $3$ is defined from a homogeneous potential $W.$ In this paper, we prove that if we add to $W$ any potential of smaller degree, we get a Poincar\'e-Birkhoff-Witt deformation of $A.$ Such PBW deformations are Calabi-Yau and are characterised among all the PBW deformations of $A.$ Various examples are presented.
\end{abstract}


\section{Introduction}
The notion of Calabi-Yau algebra has its origin in algebraic geometry: a smooth projective variety is Calabi-Yau if its canonical bundle is trivial, \textit{i.e.} if it admits a global volume form (see \textit{e.g.} \cite{nee:abel}). By Serre duality, it is equivalent to say that the bounded derived category of coherent sheaves is endowed with a Serre functor (see \cite{Bondal-Kapranov}) given by a power of the shift functor. This condition was used by Kontsevich to define the notion of a Calabi-Yau triangulated category. Given a (noncommutative) algebra $A$ over a field $\kk,$ Calabi-Yau conditions for $A$ can be defined by considering various triangulated categories associated with $A,$ for example the bounded derived category (or the stable category) of finite dimensional $A$-modules (that is, of $A$-modules that are finite dimensional as $\kk$-vector spaces). Such categories, and the noncommutative algebras giving rise to them, are currently  playing an important role in representation theory and in  algebraic geometry (resolutions of singularities) in connection with cluster algebras or McKay correspondence \cite{kelrei:clust, iyarei:mut, vg:project}.

String theory physicists have a construction of noncommutative Calabi-Yau algebras, 
which are quiver algebras with relations \cite{vdb:pot}. For a fixed quiver $Q$, this construction only depends on a noncommutative polynomial (or series) in several variables called \textit{potential}.
 Van den Bergh conjectured that any 3-dimensional Calabi-Yau algebra derives from a 
potential \cite{vdb:pot}. This was actually proved in the graded situation by  Bocklandt \cite{bock:gcy}. Our paper deals with deformations of the graded algebras $A(Q, W_{N+1})$ obtained by Bocklandt. Here  $W_{N+1}$ denotes the potential (which can be chosen  homogeneous since $A$ is graded) and $N+1$ is its degree. 

Our first observation is that the graded algebras $A(Q, W_{N+1})$ are $N$-Koszul, so we 
can consider Poincar\'e-Birkhoff-Witt (PBW) deformations of $A(Q, W_{N+1})$ in the sense of \cite{rbvg:hisym} (see also \cite{fv:pbwdef} for the case in which $Q$ has only one vertex). The main result of this paper is the following.
\Bte
Let $A=A(Q, W_{N+1})$ be a graded Calabi-Yau algebra of dimension 3. Let $W'$ be a not necessarily homogeneous potential of degree at most $ N$. Then 

(i) $A':=A(Q,W_{N+1}+W')$ is a PBW deformation of $A$, and $A'$ is Calabi-Yau. 

(ii) Assume that the characteristic of the ground field does not divide $N!$. If $A'$ is any PBW deformation of $A$, then $A'$ derives from a potential if and only if a certain condition (PBW2') (involving only the  degree $N-1$ part of the relations of $A'$) holds. 
\Ete

This theorem is the combination of Theorems  \ref{nonhomogisPBW}, \ref{NSCpbw2'} and \ref{PBWdefCY} in the text below. Part $(ii)$ can be viewed as an answer to an analogue of Van den Bergh's conjecture in a filtered situation. 

We illustrate this theorem in Section \ref{examples} by several examples arising from various sources: Yang-Mills algebras (from theoretical physics), cubic Artin-Schelter regular (AS-regular) algebras (from noncommutative algebraic geometry), antisymmetriser algebras (from representation theory), and some quiver algebras with several vertices (introduced in \cite{bock:gcy}). 

The definition of Calabi-Yau algebras used in this text is due to Ginzburg \cite{vg:project};  we relate it to that used by Bocklandt in Section \ref{generalities} (see also \cite[Section 3.2]{EO} where a similar class of algebras is introduced). Ginzburg's definition is a noncategorical definition involving Hochschild cohomology, and inspired by Van den Bergh's duality theorem. This definition is well-adapted to the determination of Calabi-Yau algebras within the class of $N$-Koszul and Artin-Schelter Gorenstein (AS-Gorenstein) algebras (Proposition \ref{auto} below). In order to explain this, we need some general considerations on graded Calabi-Yau algebras presented in Subsection \ref{grCYASG}.

We would like to discuss here the following point that appears to us to be important and is a consequence of our results: the class $\C_1$ of graded Calabi-Yau algebras of dimension $3$ is strictly contained in the class $\C_2$ of graded $N$-Koszul and AS-Gorenstein algebras of global dimension $3$, and Proposition \ref{auto} specialised to global dimension $3$ enables us to characterise $\C_1$ inside $\C_2.$ For instance, an AS-regular algebra of global dimension $3$ is in $\C_2,$ and moreover it is in $\C_1$ if and only if it is of type {\rm A} (Proposition \ref{typeofASregular}). In \cite{mdv:multi,dubois-violette}, Dubois-Violette describes the algebras inside $\C_2$ using multilinear forms called \textit{$3$-regular}, and he conjectures that any $3$-regular multilinear form provides an algebra in $\C_2.$ Therefore he offers, in a more general situation (but for quivers with only one vertex), a notion equivalent to that of a ``good''
 potential (that is, a potential that defines a Calabi-Yau algebra) with an explicit description. Moreover, he suggests an analogous description in dimension greater than $3.$

\medskip

\textbf{Acknowledgements:} We are grateful to Victor Ginzburg, Bernhard Keller, Bernard Leclerc and Michel Van den Bergh for their comments and discussions. We also thank Mariano Suarez-\'Alvarez for pointing out reference \cite{EO}.

%%%%%%%%%%%%%%%%%%%%%%%%%%%%%%%%%%%%%%%%%%%%%%%%%%%%%%%%%%%%%%%%%%%%%%%%%%%%%%%%%%%%%%%%%%%%%%%%%%

\section{Preliminaries}\label{generalities}

Let $\kk$ be a field, and let  $A$ be an  associative $\kk$-algebra which has a finite projective $A$-bimodule resolution by bimodules of finite type. The symbol $\ot$ with no subscript denotes the tensor product over the base field $\kk.$

The space $A\ot A$ is endowed with two $A$-bimodule structures: the \textit{outer} structure defined by $a\cdot(x\ot y)\cdot b=ax\ot yb,$ and the \textit{inner} structure defined by $a\cdot(x\ot y)\cdot b=xb\ot ay.$ Consequently, the Hom spaces $\Hom_{A-A}(M,A\ot A)$ of $A$-bimodule morphisms from $M$ to $A\ot A$ endowed with the outer structure are again $A$-bimodules using the inner structure of $A\ot A,$ and the same is true of the Hochschild cohomology spaces $\hh^k(A,A\ot A).$   

We use Ginzburg's definition \cite{vg:project} of a Calabi-Yau algebra:

\bd\label{ginzburgcy} We say that an algebra  $A$ as above is a \textit{Calabi-Yau algebra} of dimension $n\pgq 1$ if there are $A$-bimodule isomorphisms $\hh^k(A,A\ot A)\cong
\begin{cases}
A \mbox{ if $k=n$}\\0 \mbox{ otherwise.}
\end{cases}
$
\ed

\bp\label{hochschilddim} If $A$ is a Calabi-Yau algebra of dimension $n$, then the Hochschild dimension of $A$ (that is, the projective dimension of $A$ as an $A$-bimodule) is $n.$
\ep

\bpf Since $A$ has a finite projective $A$-bimodule resolution by assumption, the Hochschild dimension of $A$ is finite. This Hochschild dimension is at least $n,$ since $\Ext_{A-A}^n(A,A\ot A)\neq 0$ (see \cite[VI.2.1]{CE}). 

By assumption, there exists a finite projective $A$-bimodule resolution of $A$ by bimodules of finite type. Let $m$ be the shortest length of such a resolution, and let $0\rightarrow P_m\stackrel{\delta_m}{\rightarrow} P_{m-1}\stackrel{\delta_{m-1}}{\rightarrow} \cdots \rightarrow P_0\rightarrow A \rightarrow 0$ be  a projective $A$-bimodule resolution of $A$ by bimodules of finite type of length $m$. Clearly, $m$ is greater than the Hochschild dimension of $A,$ so we need only prove that $m=n.$ Assume for a contradiction that $m>n.$ Then $H^m(A,A\ot A)=0,$ and therefore $H^m(A,P)=0$ for any projective $A$-bimodule of finite type $P.$ In particular, $H^m(A,P_m)=0,$ so the map $\Hom_{A-A}(P_{m-1},P_m)\stackrel{-\circ\delta_m}{\longrightarrow} \Hom_{A-A}(P_m,P_m)$ is onto. Therefore there exists $\sigma\in \Hom_{A-A}(P_{m-1},P_m)$ such that $\sigma\circ\delta_m=\id_{P_m},$ \textit{i.e.} $\delta_m$ splits, so $P_m$ is a direct summand in $P_{m-1}. $ Therefore we can write $P_{m-1}=P_m\oplus Q$ with $Q$ a projective $A$-bimodule of finite type. But then the sequence $0\rightarrow Q\stackrel{{\delta_{m-1}}_{\mid Q}}{\longrightarrow}P_{m-2}\rightarrow\cdots \rightarrow P_0\rightarrow A\rightarrow 0$ is a projective $A$-bimodule resolution of $A$ by bimodules of finite type of length strictly less than $m,$ a contradiction. Therefore $m=n$ and the Hochschild dimension of $A$ is $n.$
\epf

We shall use the results and notations of Bocklandt \cite{bock:gcy}. We recall some of them here.

Let $Q$ be a connected quiver, fixed throughout. Let $Q_j$ denote the set of paths in $Q$ of length $j$. We assume throughout that $Q_0$ and $Q_1$ are finite. The maps $s,t:Q_1\rightarrow Q_0$ send an arrow to its \textit{source} and \textit{target}. 

A \textit{potential} is an element in the vector space $\pot{Q}:=\kk Q/[\kk Q,\kk Q].$ It can be viewed as an element in $\kk Q$ via the map $\ci: \pot{Q}\rightarrow \kk Q$ that sends a cycle $a_n\cdots a_1$ (read from right to left) to $\sum_{i}a_{i-1}\cdots a_1a_n\cdots a_i$ (this map is denoted by $\circlearrowright$ in \cite{bock:gcy}). If $p=a_n\cdots a_2a_1$ is a path and $b$ is an arrow, then $pb^{-1}$ denotes the path $a_n\cdots a_2$ if $b=a_1$, and is $0$ otherwise. Define similarly the path $b^{-1}p.$ For each $a\in Q_1,$ we shall consider the map $\partial_a:\pot{Q}\rightarrow\kk Q$ that sends an element $p$ in $\pot{Q}$ to $\ci(p)a^{-1}=a^{-1}\ci(p).$ Finally, if $W$ is a potential, we shall denote by $A(Q,W)$ the algebra $\kk Q/\I(\set{\partial_aW\,\mid\,a\in Q_1})$ (called a vacualgebra in \cite{bock:gcy}). Bocklandt proved in \cite[Theorem 3.1]{bock:gcy} that if a quiver algebra $A=\kk Q/I$ is a graded Calabi-Yau algebra of dimension $3$, then there exists a (homogeneous) potential $W$ such that $A=A(Q,W).$ Moreover, if $W$ is a potential, the algebra $A(Q,W)$ is Calabi-Yau of dimension $3$ if and only if a certain complex $C_W$ is exact \cite[Section 4.2]{bock:gcy}.

The results in \cite{bock:gcy} are still true with the definition of a Calabi-Yau algebra given here. Indeed, Bocklandt uses only \cite[Property 2.2]{bock:gcy} applied to graded algebras, which is still true with Ginzburg's definition of a Calabi-Yau algebra: the fact that the global dimension of graded Calabi-Yau algebras of dimension $n$ is $n$ follows from Remark \ref{globaldim}, and we have the  result that follows.

\bp\label{ginzburgCY} Let $A$ be a Calabi-Yau algebra of dimension $n$ over a field, and let $X$ and $Y$ be two finite-dimensional left $A$-modules. Let $(-)^*$ denote the ordinary $\kk$-dual. Then $\Ext_A^\bullet(X,Y)^*\cong\Ext_A^{n-\bullet}(Y,X).$
\ep

The proof of this Proposition relies on the  lemmas below.

\bl\label{projresol} Let $P_\bullet\rightarrow  A\rightarrow  0$ be a projective $A$-bimodule resolution of $A.$ Then for any left $A$-module $X,$ the complex $P_\bullet\ot_A X\rightarrow  X\rightarrow  0$ is a projective left  $A$-module resolution of $X.$
\el

\bpf Each module $P_k\ot_AX$ is a projective left $A$-module, since by adjunction there is an isomorphism $\Hom_A(P_k\ot_AX,\cdot)\cong\Hom_{A-A}(P_k,\Hom_\kk(X,\cdot))$, and the functor on the right is exact since $P_k$ is projective as an $A$-bimodule.

We must now prove that $P_\bullet\ot_A X\rightarrow  X\rightarrow  0$ is exact: since $\kk$ is a field, the algebra $A$ is flat over $\kk,$ and hence the enveloping algebra $A^e=A\ot A^{op}$ is flat as a left $A$-module and as a right $A$-module. Consequently, the projective $A^e$-modules $P_k$ are flat as left $A$-modules and as right $A$-modules. Flatness is a property preserved by taking kernels of epimorphisms, so that all the syzygies $\Omega^k_{A^e}(A)$ are flat as left $A$-modules and right $A$-modules. It then follows that for all $k\in\nn,$ the sequence $0\rightarrow  \Omega^{k+1}_{A^e}(A)\ot_A X\rightarrow P_k\ot_A X\rightarrow  \Omega^k_{A^e}(A)\ot_A X\rightarrow  0$ is exact, and hence the sequence $P_\bullet \ot_AX$ is exact. 
\epf

We define the following structures: if $X$ is a left $A$-module, then $X^*$ is a right $A$-module with $\alpha^a:x\mapsto \alpha(ax)$ (for $\alpha\in X^*,$ $a\in A$ and $x\in X$); if $X$ is an $A$-bimodule, then $X^*$ is an $A$-bimodule with $^a\alpha^b:x\mapsto \alpha(bxa)$ (for $\alpha\in X^*,$ $a,b\in A$ and $x\in X$); if $X$ is a left $A$-module and $X'$ is a right $A$-module, then $X\ot X'$ is an $A$-bimodule with $a\cdot (x\ot x')\cdot b=ax\ot x' b$ (for $a,b\in A$, $x\in X$ and $x'\in X'$); and if $X$ and $Y$ are left $A$-modules, $\Hom_\kk(X,Y)$ is an $A$-bimodule with $^af^b:x\mapsto af(bx)$ (for $a,b\in A,$ $f\in \Hom_\kk(X,Y)$ and $x\in X$).

\bl  Let $X$ and $Y$ be finite-dimensional  left $A$-modules. Then there is an isomorphism  $\Ext_A^k(X,Y)\cong \hh^k(A,\Hom_\kk(X,Y)).$
\el

\bpf The result is in \cite[IX.4.4]{CE}, and can also be proved using Lemma \ref{projresol} above.\epf

\bl For any $A$-bimodule $M,$ we have $(\hh_k(A,M))^*\cong \hh^k(A,M^*).$
\el

\bpf Let $P_\bullet$ be a projective $A$-bimodule resolution. Then $\hh_\bullet(A,M)$ is the homology of $P_\bullet \ot_{A^e}M\rightarrow  0.$ Dualising this gives $0\rightarrow  (P_\bullet\ot_{A^e}M)^*.$ We have $$(P_\bullet\ot_{A^e}M)^*=\Hom_\kk(P_\bullet\ot_{A^e}M,\kk)\cong \Hom_{A^e}(P_\bullet,\Hom_\kk(M,\kk))=\Hom_{A^e}(P_\bullet,M^*)$$ with adjunction. The differentials correspond via this isomorphism, therefore the cohomology of the dual complex is $\hh^\bullet(A,M^*).$
\epf

\bl Let $X$ and $Y$ be any finite-dimensional left $A$-modules. Then there is an isomorphism of $A$-bimodules $\Hom_\kk(X,Y)^*\cong \Hom_\kk(Y,X)$.
\el

\bpf We combine the following well-known $\kk$-isomorphisms, which are also $A$-bimodule homomorphisms:
\begin{enumerate}[$\bullet$]
\item $
\begin{array}{rcl}
Y\ot X^*&\rightarrow &\Hom_\kk(X,Y)\\y\ot\alpha&\mapsto &\left[x\mapsto \alpha(x)y\right]
\end{array}
$
\item $
\begin{array}{rcl}
X\ot Y^*&\rightarrow &(Y\ot X^*)^*\\
x\ot \beta &\mapsto & \left[y\ot \alpha\mapsto \alpha(x)\beta(y)\right]
\end{array}
$
\end{enumerate}
\epf

We can now prove Proposition \ref{ginzburgCY}.

\bpf We first apply Van den Bergh's duality theorem \cite{vandenbergh}: $A$ has  a finite projective resolution by $A$-bimodules of finite type; $\hh^k(A,A\ot A)=0$ for all $k\neq n$; and $U:=\hh^n(A,A\ot A)\cong A$ is an invertible $A$-bimodule (with $U^{-1}=A$). Therefore Van den Bergh's duality theorem applies, and we have: $\hh^k(A,M)\cong \hh_{n-k}(A,M).$

Applying the lemmas above, we get: 
\begin{align*}
\Ext_A^k(X,Y)^*\cong& \hh^k(A,\Hom_\kk(X,Y))^*\cong \hh_{n-k}(A,\Hom_\kk(X,Y))^*\\&\cong \hh^{n-k}(A,\Hom_\kk(X,Y)^*)\cong \hh^{n-k}(A,\Hom_\kk(Y,X))\cong \Ext_A^{n-k}(Y,X).
\end{align*}\epf

\br\label{globaldim} Let $A$ be a Calabi-Yau algebra of dimension $n$ such that there exists a nonzero finite-dimensional $A$-module $X.$ Then the global dimension of $A$ is $n:$ we know from \cite[IX.7.6]{CE} that the global dimension of $A$ is at most the Hochschild dimension of $A,$ that is, $n.$ Moreover, by Proposition \ref{ginzburgCY}, we have $\Ext_A^n(X,X)=\Ext^0_A(X,X)^*=\Hom_A(X,X)^*\neq 0,$ so that the global dimension of $A,$ which is larger than the projective dimension of the $A$-module $X$, is at least $n.$ Therefore it is equal to $n.$ 
\er

%%%%%%%%%%%%%%%%%%%%%%%%%%%%%%%%%%%%%%%%%%%%%%%%%%%%%%%%%%%%%%%%%%%%%%%%%%%%%%%%%%%%%%%%%%%%%%%%%%%%%%%%%%%%%%%%%%%%%%%%%%%%%%%%%%%%%%%%%%%%%%%%%%%%%%%%%%%%%%%%%%%%%%%%%%%%%%%%%%%%%%%%%%%%%%%%%%%%%%%%%%%%%%%%%%%%%%%%%%%%%%%%%%%%%%%%%%%%%%%%%%%%%%%%%%%%%%%%%%%%%%%%%%%%%%%%%%%%%%%%%%%%%%%%%%%%%%%%%%%%%%%%%%%%%%%%%%%%%%%%%%%%%%%%%%%%%%%%%%%%%%%%%%%%%%%%%%%%%%%%%%%%%%%%%%%%%%%%%%%%%%%%%%%%%%%%
%%%%%%%%%%%%%%%%%%%%%%%%%%%%%%%%%%%%%%%%%%%%%%%%%%%%%%%%%%%%%%%%%%%%%%%%%%%%%%%%%%%%%%%%%%%%%%%%%%

\section{PBW deformations of graded Calabi-Yau algebras of dimension 3}
\subsection{Algebras associated with nonhomogeneous potentials}

\bt\label{nonhomogisPBW}
Let $A=A(Q,W_{N+1})$ be a graded Calabi-Yau algebra of dimension $3$, where $W_{N+1}$ is a homogeneous potential of degree $N+1$, and let $W=W_{N+1}+W'=W_{N+1}+W_N+\cdots+W_0$ be a potential with $\deg W_j=j$ for each $0\ppq j\ppq N+1.$  Then $A':=A(Q,W)$ is a PBW deformation of $A.$
\et

\bpf Denote by $R$ the $\kk Q_0$-bimodule spanned by $\set{\partial_aW_{N+1}\,\mid\, a\in Q_1}$ (which can be identified as a set with the $\kk$-vector space generated by $\set{\partial_aW_{N+1}\,\mid\, a\in Q_1}$) and by $P$ the $\kk Q_0$-bimodule spanned by $\set{\partial_aW\,\mid\, a\in Q_1}$ (similarly, this can be identified with the $\kk$-vector space generated by $\set{\partial_aW\,\mid\, a\in Q_1}$). Define $\varphi:R\rightarrow   \kk Q$ by setting $\varphi(\partial_aW_{N+1})=-\partial_aW'.$ Define also $\varphi_j:R\rightarrow   (\kk Q_1)^{\ots j}=\kk Q_j$ by composing $\varphi$ with the natural projection from $\kk Q$ to $\kk Q_j$ (so that $\varphi_j(\partial_aW_{N+1})=-\partial_aW_{j+1}$). Set $\varphi_j^{1,N}=\varphi_j\ots\id:R\ots \kk Q_1\rightarrow  \kk Q_{j+1}$ and $\varphi_j^{2,N+1}=\id\ots\varphi_j:\kk Q_1\ots R\rightarrow  \kk Q_{j+1}$. Finally, set $F^j=\bigoplus_{k=0}^j \kk Q_k.$

The $\kk$-algebra $A$ can be viewed as a $\kk Q_0$-algebra, with $\kk Q_0$ semisimple. Now since $A$ is Calabi-Yau,  the complex $C_{W_{N+1}}$ in \cite{bock:gcy} is a projective $A$-bimodule resolution of $A,$ from which we can see that $A$ is $N$-Koszul (as a $\kk Q_0$-algebra).  Therefore the  results in $\cite{rbvg:hisym}$ apply, so that $A'$ is a PBW deformation of $A$ if and only the following conditions are satisfied:

\begin{enumerate}[({PBW}1)]
\item $P\cap F^{N-1}=\set{0}$
\item $(\varphi_{N-1}^{1,N}-\varphi_{N-1}^{2,N+1})(\kk Q_1\ots R\cap R\ots \kk Q_1 )\subseteq R$
\item $\forall j,$ $1\ppq j\ppq N-1,$ $$\left(\varphi_j(\varphi_{N-1}^{1,N}-\varphi_{N-1}^{2,N+1})+\varphi_{j-1}^{1,N}-\varphi_{j-1}^{2,N+1}\right)(\kk Q_1\ots R\cap R\ots \kk Q_1)=\set{0}$$
\item $\varphi_0\left(\varphi_{N-1}^{1,N}-\varphi_{N-1}^{2,N+1}\right)(\kk Q_1\ots R\cap R\ots \kk Q_1)=\set{0}.$
\end{enumerate}

Let us first give an alternative description of $R\ots \kk Q_1\cap \kk Q_1\ots R:$ consider the $\kk$-linear map $\theta:\kk Q_0\rightarrow \kk Q_1\ots R
$ that sends a vertex $e$ to $\theta(e)=\sum_{a\in eQ_1}a\partial_aW_{N+1}.$ For any cycle $p,$ we have $$\sum_{a\in eQ_1}a\partial_ap=\sum_{a\in eQ_1\mid a\in p}aa^{-1}\ci(p)=e\ci(p)=\ci(p)e=\sum_{b\in Q_1e\mid b\in p}\ci(p)b^{-1}b=\sum_{b\in Q_1e}\partial_bp\,b.$$ Therefore $\theta(e)$ is in $R\ots \kk Q_1,$ so we have a $\kk Q_0$-linear map $\theta:\kk Q_0\rightarrow  \kk Q_1\ots R\cap R\ots \kk Q_1,$ and $\theta(e)=\sum_{a\in eQ_1}a\partial_aW_{N+1}=\sum_{b\in Q_1e}\partial_bW_{N+1}\,b=e\ci(W_{N+1}).$

The map $\theta$ is one-to-one: if $\theta(\sum_{s\in Q_0}\lambda_s s)=0,$ then for a fixed $e\in Q_0$ we have $$0=e\theta(\sum_{s\in Q_0}\lambda_s s)=\sum_{s\in Q_0}\lambda_s se\,\ci(W_{N+1})=\lambda_e e\,\ci(W_{N+1})=\lambda_e \theta(e),$$ so we only need to check that $\theta(e)\neq 0.$ We use \cite[Theorem 3.1]{bock:gcy}: $A$ is Calabi-Yau, so we know that $e$ is (for instance) the source of an arrow $b$, and moreover $b$ is contained in $\ci(W_{N+1}),$ therefore $e$ is contained in $\ci(W_{N+1})$, so $\theta(e)=e\,\ci(W_{N+1})\neq 0.$

To prove that $\theta$ is an isomorphism, we shall prove that $\kk Q_0$ and $R\ots \kk Q_1\cap \kk Q_1\ots R$ have same dimension. We know from \cite{rbvg:hisym} that $R\ots \kk Q_1\cap \kk Q_1\ots R=\left(\Tor^A_3(\kk Q_0,\kk Q_0)\right)_{N+1}=\Tor^A_3(\kk Q_0,\kk Q_0)$, so $\dim_\kk R\ots \kk Q_1\cap \kk Q_1\ots R=\dim_\kk\Tor^A_3(\kk Q_0,\kk Q_0).$ As before (\textit{c.f.} Lemma \ref{projresol}), $C_{W_{N+1}}\ot_{A}\kk Q_0$ is a projective left $A$-module resolution of $\kk Q_0.$ Applying the functor $\kk Q_0\ot_A-$ gives $\Tor^A_3(\kk Q_0,\kk Q_0)=\ker(\id\ot_A \delta_3\ot_A\id)$ using the notations in $\cite{bock:gcy}.$ Therefore 
\begin{align*}
\# Q_0&\ppq \dim_\kk R\ots \kk Q_1\cap \kk Q_1\ots R \\&= \dim_\kk\Tor^A_3(\kk Q_0,\kk Q_0)\\
&\ppq \dim_\kk \kk Q_0\ot_A A\ots \mathrm{span}_\kk\set{e\ci(W_{N+1})\,\mid\, e\in Q_0}\ots A\ot_A \kk Q_0\\
&=\dim_\kk\mathrm{span}_\kk\set{e\ci(W_{N+1})\,\mid\, e\in Q_0}\ppq \# Q_0
\end{align*} so finally $\dim_\kk R\ots \kk Q_1\cap \kk Q_1\ots R=\dim_\kk \kk Q_0$ and $\theta$ is an isomorphism. In particular, $\set{\theta(e)\,\mid\,e\in Q_0}$ is a basis for $R\ots \kk Q_1\cap \kk Q_1\ots R.$

Let us now check that conditions (PBW1) to (PBW4) hold:
\begin{enumerate}[({PBW}1)]
\item Clearly $P\cap F^{N-1}=\set{\sum_{a\in Q_1}\lambda_a \partial_aW\,\mid\, \lambda_a\in\kk \mbox{ and } \sum_{a\in Q_1}\lambda_a \partial_aW_{N+1}=0}.$ However, the $\partial_aW_{N+1}$ for $a\in Q_1$ are linearly independent: consider the resolution $C_{W_{N+1}}$ in \cite{bock:gcy}. By Lemma \ref{projresol}, $C_{W_{N+1}}\ot_A \kk Q_0$ is a projective resolution of $\kk Q_0$ as a left $A$-module. The maps in the complex $\Hom_A(C_{W_{N+1}}\ot_A \kk Q_0,\kk Q_0)$ vanish, so $\# Q_1=\dim_\kk\Ext_A^1(\kk Q_0,\kk Q_0)=\dim_\kk\Ext_A^2(\kk Q_0,\kk Q_0)=\dim_\kk\mathrm{span}_\kk\set{\partial_aW_{N+1}\,\mid\, a\in Q_1}$ using the Calabi-Yau property. Therefore, $P\cap F^{N-1}=\set{0}.$

\item This condition is equivalent to: $\forall e\in Q_0,$ $(\varphi_{N-1}^{1,N}-\varphi_{N-1}^{2,N+1})(\theta(e))\in R.$ We have 
\begin{align*}
(\varphi_{N-1}^{1,N}-\varphi_{N-1}^{2,N+1})(\theta(e))&=\varphi_{N-1}^{1,N}\left(\sum_{b\in Q_1 e}\partial_bW_{N+1}\,b\right)-\varphi_{N-1}^{2,N+1}\left(\sum_{a\in eQ_1}a\partial_aW_{N+1}\right)\\&
=\sum_{b\in Q_1 e}\varphi_{N-1}(\partial_bW_{N+1})b-\sum_{a\in eQ_1}a\varphi_{N-1}(\partial_aW_{N+1})\\
&=-\sum_{b\in Q_1 e}\partial_bW_{N}\,b+\sum_{a\in eQ_1}a\partial_aW_{N}=0\in R
\end{align*} as above.

\item Using the fact that $(\varphi_{N-1}^{1,N}-\varphi_{N-1}^{2,N+1})(\theta(e))=0$ for all $e\in Q_0$, this condition becomes $(\varphi_{j-1}^{1,N}-\varphi_{j-1}^{2,N+1})(\theta(e))=0$ for all $e\in Q_0.$ This is true by a  calculation similar to that above.

\item This condition is trivially true using the fact that $(\varphi_{N-1}^{1,N}-\varphi_{N-1}^{2,N+1})(\theta(e))=0$ for all $e\in Q_0$.
\end{enumerate} 

Therefore $A'$ is a PBW deformation of $A.$
\epf

This leads to the following question: when is a PBW deformation of a Calabi-Yau algebra of dimension $3$ obtained from a (nonhomogeneous) potential? 

Let $A=A(Q,W_{N+1})$ be a Calabi-Yau algebra of dimension $3,$ and let $A'$ be a PBW deformation of $A$. We recall here the notations of \cite{rbvg:hisym}: let $P$ be a sub-$\kk Q_0$-bimodule of $F^N$ such that $A'=\kk Q/\I(P),$ let $\pi:F^N\rightarrow\kk Q_N$ be the natural projection, and set $R=\pi(P),$ with $A=\kk Q/\I(R)$. The condition $P\cap F^{N-1}=\set{0}$ holds since $A'$ is a PBW deformation of $A,$ so $\pi:P\rightarrow R$ is an isomorphism. For each $a\in Q_1,$ there exists a unique element in $P$ whose image by $\pi$ is $\partial_aW_{N+1}$; this element can be written $\partial_aW_{N+1}-\varphi(\partial_aW_{N+1})$ with $\varphi(\partial_aW_{N+1})\in F^{N-1}$. This defines a $\kk Q_0$-linear map $\varphi:R\rightarrow F^{N-1}, $ and $P=\mathrm{span}_{\kk Q_0}\set{\partial_aW_{N+1}-\varphi(\partial_aW_{N+1})\,\mid\,a\in Q_1}.$ Define maps $\varphi_j:R\rightarrow kQ_j$ as above by composing $\varphi$ with the projections $\kk Q\rightarrow\kk Q_j$ for $0\ppq j\ppq N-1.$

 If $A'$ can be obtained from a potential, that is, if there exists a potential $W=W_{N+1}+W'$ with $\deg W'\ppq N$ such that $A'=A(Q,W),$ then we have $\varphi(\partial_aW_{N+1})=-\partial_aW',$ and therefore  $$\mbox{(PBW2') }\ \forall e\in Q_0,\  (\varphi_{N-1}^{1,N}-\varphi_{N-1}^{2,N+1})(\theta(e))=0,$$ since $(\varphi_{N-1}^{1,N}-\varphi_{N-1}^{2,N+1})(\theta(e))=-\sum_{b\in Q_1 e}\partial_bW_{N}\,b+\sum_{a\in eQ_1}a\partial_aW_{N}=0.$ This condition (PBW2') is in fact sufficient:

\bt\label{NSCpbw2'} Let $A=A(Q,W_{N+1})$ be a Calabi-Yau algebra of dimension $3,$ and let $A'$ be a PBW deformation of $A$. Assume that $\car(\kk)$ does not divide $N!$. Then  $A'$ can be obtained from a potential if and only if condition (PBW2') holds.
\et

\bpf 

Let $A'$ be a PBW deformation of $A$ satisfying (PBW2'). Since $\varphi_{j-1}$ is $\kk Q_0$-linear,  $\varphi_{j-1}(\partial_aW_{N+1})$ is a linear combination of paths all starting at $t(a)$ and all ending at $s(a).$ Therefore $A'$ satisfies,  for all $1\ppq j\ppq N$, $$(*_j)\ \forall e\in Q_0, \ \sum_{a\in eQ_1}a\varphi_{j-1}(\partial_aW_{N+1})-\sum_{b\in Q_1e}\varphi_{j-1}(\partial_bW_{N+1})b=0.$$

The element $\varphi_{j-1}(\partial_aW_{N+1})$ is a linear combination of paths of length $j-1$ starting at $t(a)$ and ending at $s(a)$, therefore we can write $$\varphi_{j-1}(\partial_aW_{N+1})=\sum_{q\in s(a)Q_{j-1}t(a)}\hspace*{-.3cm}\lambda_{a,q}q=\sum_{q\in Q_{j-1}\mid aq\in t(a)Q_jt(a)}\hspace*{-.6cm}\lambda_{a,q}q=\sum_{q\in Q_{j-1}\mid qa\in s(a)Q_js(a)}\hspace*{-.6cm}\lambda_{a,q}q$$ for some scalars $\lambda_{a,q}\in \kk.$  Then $(*_j)$ becomes $$\forall e\in Q_0, \ \sum_{\tiny
\begin{array}{c}
a\in eQ_1,\\ q\in Q_{j-1}\mid aq\in eQ_je
\end{array}
}\hspace*{-.6cm}\lambda_{a,q}aq=\sum_{\tiny
\begin{array}{c}
b\in Q_1e,\\q\in  Q_{j-1}\mid qb\in eQ_je
\end{array}
 }\hspace*{-.6cm}\lambda_{b,q}qb.$$

If $\sigma$ is a cycle of length $j\pgq 2,$ then we can write $\sigma=aqb$ with $a,b\in Q_1 $ and $q\in Q_{j-2}.$ This cycle occurs once on each side of the equality above with $e=t(a)=s(a)$, so $\lambda_{a,qb}=\lambda_{b,aq}$. Moreover, if $j=1,$ we also have $\lambda_{a,e}=\lambda_{a,e}$ if $a\in eQ_1e.$

Now if we write $\sigma=a_j\ldots a_1,$ we have $\lambda_{a_j,a_{j-1}\ldots a_1}=\lambda_{a_1,a_j\ldots a_2}=\lambda_{a_2,a_1a_j\ldots a_3}=\cdots,$ so that the coefficient only depends on $\ci(\sigma).$ We denote this common coefficient by $\hat{\lambda}_\sigma.$

Let $\sigma$ be a cycle in $Q.$ Define a new element $\ci'(\sigma)$ in $\kk Q$ as follows: if $\sigma$ is not a power of a cycle, then we set $\ci'(\sigma)=\ci(\sigma);$ if $\sigma=\tau^m$ is a power of a cycle, with $m$ maximal, then $\ci'(\sigma)$ is the sum of the $m^{\tiny\rm th}$ powers of all the terms in $\ci(\tau).$ We shall denote by $\bar{\sigma}$ the class of a cycle $\sigma$ modulo cyclic permutations.

Now define $W_j:=-\frac{1}{j}\sum_{e\in Q_0}\sum_{\sigma\in eQ_je}\hat{\lambda}_\sigma\ci'(\sigma).$ Then, if $a\in Q_1,$ we want to compute $\partial_a W_j.$ For this, we shall need the following lemmas.

\bl\label{technicallemma} Let $\sigma$ be a cycle of length $j$ which can be written $\sigma=ap$ where $a$ is an arrow and $p$ is a path, and such that there exists a cyclic permutation of length at most $j-1$ which when applied to $\sigma$ again gives $ap.$ Then $\sigma$ is a power of a cycle.
\el

\bpf The hypothesis implies that $a$ occurs in $p,$ so we can write $p=\alpha a\beta$ for some paths $\alpha$ and $\beta,$ and these paths can be chosen such that $\sigma=a\alpha a\beta=a\beta a\alpha\ (\circledast).$

Let $\ell(q)$ denote the length of a path $q.$ If $\ell(\beta)=\ell(\alpha),$ then necessarily $\alpha=\beta$ and therefore $\sigma=(a\alpha)^2$ is a power of a cycle.

If $\ell(\beta)>\ell(\alpha),$ then the expressions of $\sigma$ above show that $\beta=\alpha a\beta'$ for some path $\beta',$ with $\ell(\beta')<\ell(\beta).$ Therefore, $\sigma=a\alpha a\alpha a\beta'=a\alpha a\beta' a\alpha$ and so $a\alpha a\beta'=a\beta'a\alpha.$ We have obtained an identity similar to $(\circledast),$ but with a path of smaller length. A descending induction shows that $\sigma$ is a power of a cycle.
\epf

\bl Let $\sigma$ be a cycle of length $j$ and let $a$ be an arrow. Then $\partial_a \ci'(\sigma)=j\sum_{q\in Q_{j-1}\mid\,\ovl{aq}=\bar{\sigma}}q.$
\el

\bpf 
First assume that $\sigma$ is not a power of a cycle. Then $\ci'(\sigma)=\ci(\sigma).$ This is a sum of $j$ terms, and all the terms are obtained by applying cyclic permutations to $\sigma$, and therefore they all give the same derivatives. So we need only look at the derivative with respect to $a$ of $\sigma$, and multiply the result by $j.$ 

No cyclic permutation of length at most $j-1$ of $\sigma$ is equal to $\sigma,$ and therefore if $\bar{\sigma}=\ovl{ap}=\ovl{aq}$ for two different choices of $a$ within $\sigma,$ then $p\neq q$ (Lemma \ref{technicallemma}). Therefore $\partial_a\sigma=\sum_{q\in Q_{j-1}\mid\, \ovl{aq}=\bar{\sigma}}q$ and we get $\partial_a\ci'(\sigma)=j\sum_{q\in Q_{j-1}\mid\, \ovl{aq}=\bar{\sigma}}q.$

Now assume that $\sigma$ is a power of a cycle, and write $\sigma=\tau^m$ with $m$ maximal (so that $\tau$ is not a power of a cycle). As before, there are $\frac{j}{m}$ terms in $\ci'(\sigma),$ all obtained by applying cyclic permutations to $\sigma,$ so we need only consider $\partial_a\sigma$ and multiply the result by $\frac{j}{m}.$

To compute $\partial_a\sigma=\partial_a\tau^m,$ we choose a factor $\tau,$ then for each $a$ in this factor, we apply the appropriate cyclic permutation to $\sigma$ so that the $a$ occurs at the end of the cycle, and then we remove the $a.$ We then do the same for each of the other factors $\tau,$ and in each case we will obtain the same result. Therefore a factor $m$ appears when we compute this derivative: $$\partial_a\tau^m=m\sum_{(p_1,p_2)\mid\, \tau=p_2ap_1}  p_1\tau^{m-1}p_2$$ (in this sum, $p_1p_2$ is necessarily a path of length $\frac{j}{m}-1$). Now consider the map $\varphi$ from $\set{(p_1,p_2)\,\mid\, p_2ap_1=\tau}$ to $\set{q\,\mid\, \ovl{aq}=\bar{\sigma}}$ given by $\varphi(p_1,p_2)= p_1\tau^{m-1}p_2.$  This map is one-to-one: if $\varphi(p_1,p_2)=\varphi(p_1',p_2')$ with $(p_1,p_2)\neq(p_1',p_2')$, assume that $\ell(p_1)<\ell(p_1')$, so that $\ell(p_2')<\ell(p_2).$ Since $p_1\tau^{m-1}p_2=p_1'\tau^{m-1}p_2',$ we have $p_1'=p_1p_1''$ and $p_2=p_2''p_2'$ for some paths $p_1''$ and $p_2''$, so we get $\tau^{m-1}p_2''=p_1''\tau^{m-1}.$ By \cite[Lemma 2.4]{GS}, there exist a cycle $\gamma$ and an integer $t\pgq2$ such that $\tau=\gamma^t$, a contradiction. Therefore $\varphi$ is one-to-one. The map $\varphi$ is also onto: if $\bar{\sigma}=\ovl{aq},$ the $a$ in this equality is chosen within $\sigma$ and therefore within one of the factors $\tau,$ so we can write $\tau=p_2ap_1$ and $\sigma=\tau^up_2ap_1\tau^v$ for some integers $u$ and $v$ such that $u+v=m-1,$ and therefore $q=p_1\tau^{m-1}p_2=\varphi(p_1,p_2).$

Finally, we have $\partial_a\sigma=m\sum_{q\in Q_{j-1}\mid\, \ovl{aq}=\bar{\sigma}}q,$ and this gives the result.
\epf

Now we have

\begin{align*}\partial_aW_j=-\frac{1}{j}\sum_{e\in Q_0}\sum_{\sigma\in eQ_je}\hat{\lambda}_\sigma\partial_a\ci'(\sigma)&=-\frac{1}{j}\sum_{\tiny
\begin{array}{c}
e\in Q_0,\\\sigma\in eQ_je
\end{array}}
\sum_{q\in Q_{j-1}\mid \ovl{aq}=\bar{\sigma}}\hspace*{-.5cm}j\hat{\lambda}_\sigma q\\&=-\sum_{\tiny
\begin{array}{c}
q\in Q_{j-1}\mid \\aq\in t(a)Q_jt(a)
\end{array}
}\hspace*{-.9cm}{\lambda}_{a,q} q=-\varphi_{j-1}(\partial_aW_{N+1}).\end{align*} Therefore $A'=A(Q,W)$ with $W=W_{N+1}+W_N+\cdots+W_1.$\epf

\br In the result above, the potential $W=W_{N+1}+W_N+\cdots+W_1$ defining $A'$ is unique. Indeed, it is easy to see that $\cap_{a\in Q_1}\ker(\partial_a)=\kk Q_0$, and since the degree of each $W_j$ is nonzero, it follows that two such potentials must be equal.
\er

%%%%%%%%%%%%%%%%%%%%%%%%%%%%%%%%%%%%%%%%%%%%%%%%%%%%%%%%%%%%%%%%%%%%%%%%%%%%%%%%%%%%%%%%%%%%%%%%%%%%%%%%%%%%%%%%%%%%%%%%%%%%%%%%%%%%%%%%%%%5

\subsection{The algebras $A(Q,W)$ are Calabi-Yau }

We now prove the following result.

\bt\label{PBWdefCY} Let $A=A(Q,W_{N+1})$ be a graded Calabi-Yau algebra of dimension $3$. Let $A'=A(Q,W)$ be a PBW deformation of $A$ defined by a potential $W=W_{N+1}+W'$ with $\deg W'\ppq N.$ Then $A'$ is Calabi-Yau of dimension 3.
\et

The proof will use the following characterisation of filtered Calabi-Yau algebras (the analogue of \cite[Theorem 4.2]{bock:gcy}):

\bl\label{selfdualCY} Let $B$ be a filtered $\kk$-algebra. If $B$ has a selfdual (with respect to the functor $(\cdot)\spcheck=\Hom_{B-B}(\cdot,B\ot B)$) projective $B$-bimodule resolution of finite length $n$ by $B$-bimodules of finite type,
then $B$ is Calabi-Yau of dimension $n.$
\el

\bpf By assumption, we have a diagram
$$\xymatrix{P_n\ar[d]^{\alpha_{n}}\ar[r]^{d_n}&P_{n-1}\ar[d]^{\alpha_{n-1}}\ar[r]^{d_{n-1}}&\cdots&\ar[r]^{d_{2}}&P_1\ar[d]^{\alpha_{1}}\ar[r]^{d_1}&P_0\ar[d]^{\alpha_{0}}\\
P_0\spcheck\ar[r]^{-d_1\spcheck}&P_{1}\spcheck\ar[r]^{-d_{2}\spcheck}&\cdots&\ar[r]^(.3){-d_{n-1}\spcheck}&P_{n-1}\spcheck\ar[r]^{-d_n\spcheck}&P_n\spcheck
}$$  in which all the squares commute and where all the maps $\alpha_i$ are isomorphisms.
Then:
\begin{enumerate}[$\bullet$]
\item $B$ has a finite projective $B$-bimodule resolution by bimodules of finite type.
\item For $0<k<n$, we have 
\begin{align*}
\hh^k(B,B\ot B)=\frac{\ker(P_k\spcheck\stackrel{d_{k+1}\spcheck}{\longrightarrow}P_{k+1}\spcheck)}{\im (P_{k-1}\spcheck\stackrel{d_{k}\spcheck}{\longrightarrow}P_k\spcheck)}\cong \frac{\ker(P_{n-k}\stackrel{-d_{n-k}}{\longrightarrow}P_{n-k-1})}{\im(P_{n-k+1}\stackrel{-d_{n-k+1}}{\longrightarrow}P_{n-k})}=0
\end{align*}
\item For $k>n$, $\hh^k(B,B\ot B)=0.$
\item $
\hh^n(B,B\ot B)=\dsty\frac{\ker(P_n\spcheck{\longrightarrow} 0)}{\im (P_{n-1}\spcheck\stackrel{d_{n-1}\spcheck}{\longrightarrow}P_n\spcheck)}\cong \frac{P_0}{\im(P_1\stackrel{-d_1}{\longrightarrow} P_0)}\cong B.
$
\end{enumerate} Therefore $B$ is Calabi-Yau of dimension $n$.
\epf

We now prove Theorem \ref{PBWdefCY}.

\bpf We shall construct a selfdual projective resolution of $A'.$

Let $R$ denote the $\kk Q_0$-bimodule spanned by $\set{\partial_aW_{N+1}\,\mid\,a\in Q_1}$ and $P$ the $\kk Q_0$-bimodule  spanned by $\set{\partial_aW\,\mid\,a\in Q_1}.$ Let $\pi:P\rightarrow R$ be the projection defined by $\partial_aW\mapsto\partial_aW_{N+1}.$
Consider the $\kk Q_0$-bimodule map $\tilde{\theta}:\kk Q_0\rightarrow P\ots \kk Q_1$ that sends a vertex $e$ to $\sum_{a\in Q_1e}\partial_aW\, a.$
 For every $e\in Q_0,$ $\tilde{\theta}(e)\in \kk Q_1\ots P:$ 
\begin{align*}
\tilde{\theta}(e)&=\theta(e)-\sum_{a\in Q_1e}\varphi(\partial_aW_{N+1})a=\theta(e)+\sum_{a\in Q_1e}\partial_aW'\,a\\
&=\theta(e)+\sum_{b\in eQ_1}b\partial_bW'=\theta(e)-\sum_{b\in eQ_1}b\varphi(\partial_bW_{N+1})=\sum_{b\in eQ_1}b\partial_bW.
\end{align*} Therefore $\tilde{\theta}$ takes values in $P\ots \kk Q_1 \cap \kk Q_1\ots P.$ We wish to compose $\tilde{\theta}$ with $\pi\ots\id$ and with $\id\ots\pi:$ take an element $x=\sum_{a,b\in Q_1}a\partial_bW=\sum_{c,d\in Q_1}\partial_cW\,d $ in $ P\ots \kk Q_1 \cap \kk Q_1\ots P;$ then looking at the terms of degree $N+1$ in $x$ gives $\sum_{a,b\in Q_1}a\partial_bW_{N+1}=\sum_{c,d\in Q_1}\partial_cW_{N+1}\,d,$ that is, $(\id\ots\pi)(x)=(\pi\ots \id)(x).$ Therefore $(\pi\ots\id)\tilde{\theta}=(\id\ots\pi)\tilde{\theta}$ takes values in $R\ots \kk Q_1 \cap \kk Q_1\ots R$ and is equal to $\theta,$ which is an isomorphism. We deduce that $\tilde{\theta}:\kk Q_0\rightarrow \im\tilde{\theta}$ is an isomorphism, and that $(\pi\ots\id)_{\mid \im\tilde{\theta}}$ is an isomorphism from $\im\tilde{\theta}$ to $\im\theta.$

The sequence $\C$ below is a complex of projective $A'$-bimodules of finite type: 
\begin{align*}
\C:\ 0\rightarrow A'\ots \im\tilde{\theta}\ots A'&\stackrel{\delta_3}{\rightarrow}A'\ots P\ots A'\stackrel{\delta_2}{\rightarrow}A'\ots \kk Q_1\ots A'\\&\stackrel{\delta_1}{\rightarrow}A'\ots \kk Q_0\ots A'\stackrel{\mu}{\rightarrow}A'\rightarrow 0
\end{align*} where $\mu$ is the multiplication map and  $\delta_1,\delta_2,\delta_3$ are morphisms of $A'$-bimodules defined by $$
\begin{array}{l}
\delta_1(1\ots a\ots 1)={a}\ots s(a)\ots 1-1\ots t(a)\ots {a},\\
\delta_2(1\ots \partial_aW\ots 1)=\Delta(\partial_aW), \\
\delta_3(1\ots \tilde{\theta}(e)\ots 1)=\sum_{a\in eQ_1}{a}\ots \partial_aW\ots 1-\sum_{b\in Q_1 e}1\ots \partial_bW\ots {b}.
\end{array}
$$ For $\delta_2,$ the map $\Delta$ is the linear map  defined on the monomials by $\Delta(a_k\ldots a_1)=\sum_{i=1}^k{a_k\cdots a_{i+1}}\ots  a_i\ots  {a_{i-1}\cdots a_1}$, and $\Delta(1_\kk)=0$ (the bar denotes the class of an element in $\kk Q$ modulo $P$).

Clearly, $\mu\delta_1=0.$ We have $\delta_1\circ\delta_2(1\ots\partial_aW\ots1)={\partial_aW}\ots t(a)\ots 1+\mbox{terms that simplify}-1\ots s(a)\ots {\partial_aW}=0.$ Finally, 
\begin{align*}
\delta_2\circ\delta_3(1\ots \tilde{\theta}(e)\ots 1)&=\sum_{a\in eQ_1}{a}\Delta(\partial_aW)-\sum_{b\in Q_1e}\Delta(\partial_bW){b}\\&=\left[\Delta(\tilde{\theta}(e))-\sum_{a\in eQ_1}1\ots a \ots {\partial_aW}\right]\\&\ \ \ \ \ -\left[\Delta(\tilde{\theta}(e))-\sum_{b\in Q_1e}
{\partial_bW}\ots b \ots 1\right]=0.
\end{align*}

We will now prove that $\C$ is a resolution of $A':$ the complex $\C$ is isomorphic to the complex $\C'$ defined as follows:
\begin{align*}
\C':\ 0\rightarrow A'\ots \im{\theta}\ots A'&\stackrel{\delta'_3}{\rightarrow}A'\ots R\ots A'\stackrel{\delta'_2}{\rightarrow}A'\ots \kk Q_1\ots A'\\&\stackrel{\delta_1}{\rightarrow}A'\ots \kk Q_0\ots A'\stackrel{\mu}{\rightarrow}A'\rightarrow 0
\end{align*} with $\delta'_2(1\ots \partial_aW_{N+1}\ots 1)=\Delta(\partial_aW)$ and $\delta'_3(1\ots \theta(e)\ots 1)=\sum_{a\in eQ_1}{a}\ots \partial_aW_{N+1}\ots 1-\sum_{b\in Q_1 e}1\ots \partial_bW_{N+1}\ots {b}.$ The isomorphism between the two complexes is given by the maps $\id_{A'}\ots (\pi\ots \id)\ots \id_{A'}$, $\id_{A'}\ots \pi\ots \id_{A'}$ for the first two and by $\id$ for the other three.

The complex $\C'$ is filtered, and the graded complex associated to $\C'$ is the complex $C_{W_{N+1}}$ in \cite{bock:gcy}, which is exact (since $A$ is Calabi-Yau). Therefore $\C'$ is exact (since the natural functor which associates a graded object to a filtered object is faithful), hence $\C$ is exact.

Using Lemma \ref{selfdualCY}, we need only check that $\C$ is selfdual to prove that $A'$ is Calabi-Yau of dimension $3.$ We must first define the isomorphisms $\alpha_i.$

We use Bocklandt's notation \cite{bock:gcy}: if $T$ is a finite-dimensional $\kk Q_0$-bimodule, $F_T$ is the $A'$-bimodule $A'\ots T\ots A'.$ If $T^*$ is the $\kk$-dual of $T,$ we have an isomorphism of $A'$-bimodules $F_{T^*}\stackrel{\sim}{\rightarrow}F_T\spcheck$ given by $1\ots \alpha \ots 1\mapsto \left[1\ots t \ots 1\mapsto \sum_{i,j\in Q_0}\alpha(itj)i\ot j\right]$ for $t\in T$ and $\alpha\in T^*.$

For $\alpha_0:F_{\kk Q_0}\rightarrow F_{\im\tilde{\theta}}\spcheck$, consider the composition of $\theta:\kk Q_0\rightarrow \im\tilde{\theta}$ with the $\kk$-linear isomorphism $\im\tilde{\theta}\rightarrow (\im\tilde{\theta})^*$ defined on the basis $\set{\tilde{\theta}(e)\,\mid\, e\in Q_0}$ by sending the basis elements to the corresponding elements in the dual basis. This composition is an isomorphism of $\kk Q_0$-bimodules. Now since $\kk Q_0$ is semisimple, $A'$ is flat over $\kk Q_0,$ so the functor $\id_{A'}\ots\,\cdot\,\ots\id_{A'}$ is an isomorphism of $A'$-bimodules from $F_{\kk Q_0}$ to $F_{(\im\tilde{\theta})^*}$. Composing with the isomorphism $F_{(\im\tilde{\theta})^*}\cong F_{\im{\tilde{\theta}}}\spcheck$ above gives an $A'$-bimodule isomorphism $\alpha_0:F_{\kk Q_0}\rightarrow F_{\im{\tilde{\theta}}}\spcheck.$ This isomorphism can be expressed explicitly: $$\alpha_0(1\ots e\ots 1)(1\ots \tilde{\theta}(s)\ots 1)=es\ot es.$$

For $\alpha_1:F_{\kk Q_1}\rightarrow F_{P}\spcheck,$ we proceed similarly, using the composition of the $\kk$-linear isomorphism $\kk Q_1\rightarrow P$ sending $a$ to $\partial_aW$ with the $\kk$-linear isomorphism $P\rightarrow P^*$ sending the basis $\set{\partial_aW\,\mid\,a\in Q_1}$ to its dual basis (this composition is an isomorphism of $\kk Q_0$-bimodules). We get $$\alpha_1(1\ots a\ots 1)(1\ots \partial_bW\ots 1)=\begin{cases}
s(a)\ot t(a)\mbox{ if $a=b$}\\0\mbox{ otherwise}
\end{cases}.$$

For $\alpha_2:F_{P} \rightarrow F_{\kk Q_1}\spcheck,$ we use the same construction starting with the $\kk Q_0$-bilinear map sending $\partial_aW$ to the opposite of the corresponding element in the dual basis of $\set{a\,\mid\,a\in Q_1}$ of $\kk Q_1.$ This gives $$\alpha_2(1\ots \partial_aW \ots 1)(1\ots b\ots 1)=\begin{cases}
-t(a)\ot s(a) \mbox{ if $a=b$}\\0\mbox{ otherwise}
\end{cases}.$$

Finally, the same procedure gives $\alpha_3:F_{\im\tilde{\theta}}\rightarrow F_{\kk Q_0}\spcheck,$ with $$\alpha_3(1\ots\tilde{\theta}(e)\ots1)(1\ots s\ots 1)=-es\ot es.$$

 So we have a diagram:
$$\xymatrix@=1cm{F_{\im\tilde{\theta}}\ar[d]^{\alpha_{3}}\ar[r]^{\delta_3}&F_{ P}\ar[d]^{\alpha_{2}}\ar[r]^{\delta_{2}}&F_{\kk Q_1 }\ar[d]^{\alpha_{1}}\ar[r]^{\delta_1}&F_{\kk Q_0 }\ar[d]^{\alpha_{0}}\\
(F_{ \kk Q_0})\spcheck\ar[r]^{-\delta_1\spcheck}&(F_{ \kk Q_1})\spcheck\ar[r]^{-\delta_{2}\spcheck}&(F_{ P})\spcheck\ar[r]^{-\delta_3\spcheck}&(F_{\im\tilde{\theta} })\spcheck
}$$ We must prove that it commutes.

Recall that $M\spcheck$ is the set of $A'$-bimodule morphisms from $M$ to $A'\ot A'$  when $A'\ot A'$ is endowed with the \textit{outer} bimodule structure. The result $M\spcheck$ is then an $A'$-bimodule when we endow  $A'\ot A'$  with the \textit{inner} bimodule structure.

Let $\tau:A'\ot A'\rightarrow A'\ot A'$ be the map that sends $x\ot y$ to $y\ot x.$

\begin{enumerate}[$\bullet$]
\item On the one hand,
\begin{align*}
\alpha_0&\circ\delta_1(1\ots a \ots 1)(1\ots \tilde{\theta}(e)\ots 1)\\&=\alpha_0(a\ots s(a)\ots 1-1\ots t(a)\ots a)(1\ots \tilde{\theta}(e)\ots 1)\\
&=s(a)e\ot as(a)e-t(a)ea\ot t(a)e\\
&=
\begin{cases}
e\ot a&\mbox{if $s(a)=e$}\\ -a\ot e &\mbox{if $t(a)=e$}\\0&\mbox{otherwise}
\end{cases}
\end{align*} (here it is the inner action on $A'\ot A'$ that is involved).

On the other hand, 
\begin{align*}
\delta_3\spcheck&\circ\alpha_1(1\ots a\ots 1)(1\ots \tilde{\theta}(e)\ots 1)
\\&=\alpha_1(1\ots a\ots 1)\left(\sum_{b\in eQ_1}b\ots \partial_bW\ots 1-\sum_{c\in Q_1e}1\ots \partial_cW\ots c\right)\\
&=
\begin{cases}
\alpha_1(1\ots a \ots 1)(-1\ots \partial_aW\ots a)=-e\ot a &\mbox{if $s(a)=e$}\\
\alpha_1(1\ots a \ots 1)(a\ots \partial_aW\ots 1)=a\ot e &\mbox{if $t(a)=e$}\\
0&\mbox{otherwise}
\end{cases}
\end{align*} (here it is the outer action on $A'\ot A'$ that is involved).

Therefore $\alpha_0\circ\delta_1=-\delta_3\spcheck\circ\alpha_1.$

\item Similarly, we can prove that 
 $\alpha_2\circ\delta_3=-\delta_1\spcheck\circ\alpha_3.$

\item Finally, we must prove that $\alpha_1\circ \delta_2=-\delta_2\spcheck \circ \alpha_2.$ 

Following Van den Bergh in \cite{vdb:pot} (for the one-vertex case), we define the partial derivative of a path $p$ with respect to an arrow $a$ by $\der{a}p=\sum_{p=uav}u\ot v.$ This partial derivative extends naturally to a $\kk Q_0$-linear map $\der{a}:A'\rightarrow A'\ot A'.$

Now if $p$ is a cycle and if $a$ and $b$ are arrows, we have the following relation: $$\der{a}(\partial_bp)=\tau\left[\der{b}(\partial_ap)\right].$$ Indeed,
\begin{align*}
\der{a}(\partial_bp)
&=\sum_{p=ubv}\der{a}(vu)\\
&=\sum_{p=ubv_1av_2}v_1\ot v_2u+\sum_{p=u_1au_2bv}vu_1\ot u_2\\
&=\tau\left[\sum_{p=ubv_1av_2}v_2u\ot v_1 +\sum_{p=u_1au_2bv} u_2\ot vu_1\right]\\
&=\tau\left[\sum_{p=xay}\left(\sum_{y=u_2bv}u_2\ot vx+\sum_{x=ubv_1}yu\ot v_1\right)\right]\\
&=\tau\left[\sum_{p=xay}\der{b}(yx)\right]
=\tau\left[\der{b}(\partial_ap)\right].
\end{align*} This relation extends to any element in $\pot{Q}$.

Now consider $\delta_2\spcheck\circ\alpha_2(1\ots \partial_aW\ots 1)(1\ots \partial_bW\ots 1)=\alpha_2(1\ots\partial_aW\ots 1)(\Delta(\partial_bW)).$ The element $\alpha_2(1\ots \partial_aW\ots 1)(\lambda\ots c\ots \mu)$ is nonzero only if $c=a,$ and $\alpha_2(1\ots \partial_aW\ots 1)(\lambda\ots a\ots \mu)=-\lambda\ot \mu,$ so applying $-\alpha_2(1\ots \partial_aW\ots 1)$ to $\Delta(\partial_bW)$ is the same as replacing (one at a time) the copies of $a$ inside $\partial_bW$ by the symbol $\ot$, that is, applying $\der{a}$ to $\partial_bW.$ Therefore $\delta_2\spcheck\circ\alpha_2(1\ots \partial_aW\ots 1)(1\ots \partial_bW\ots 1)=-\der{a}(\partial_bW).$

On the other hand, $\alpha_1\circ\delta_2(1\ots\partial_aW\ots 1)(1\ots \partial_bW\ots 1)=\alpha_1(\Delta\partial_aW)(1\ots \partial_bW\ots 1).$ The element $\alpha_1(\lambda\ots c\ots \mu)(1\ots \partial_bW\ots 1)$ is nonzero only if $c=b,$ and $\alpha_1(\lambda\ots b\ots \mu)(1\ots \partial_bW\ots 1)=\mu\ot \lambda=\tau(\lambda\ot \mu)$ (we must use the inner action here). So, as above, $\alpha_1\circ\delta_2(1\ots\partial_aW\ots 1)(1\ots \partial_bW\ots 1)=\tau\left[\der{b}(\partial_aW)\right]=\der{a}(\partial_bW)$.
\end{enumerate}\epf

%%%%%%%%%%%%%%%%%%%%%%%%%%%%%%%%%%%%%%%%%%%%%%%%%%%%%%%%%%%%%%%%%%%%%%%%%%%%%%%%%%%%%%%%%%%%%%%%%%%%%%%%%%%%%%%%%%%%%%%%%%%%%%%%%%%%%%%%%%%%%%%%%%%%%%%%%%%%%%%%%%%%%%%%%%%%%%%%%%%%%%%%%%%%%%%%%%%%%%%%%%%%%%%%%%%%%%%%%%%%%%%%%%%%%%%%%%%%%%%%%%%%%%%%%%%%%%%%%%%%%%%%%%%%%%%%%%%%%%%%%%%%%%%%%%%%%%%%%%%%%%%%%%%%%%%%%%%%%%%%%%%%%%%%%%%%%%%%%%%%%%%%%%%%%%%%%%%%%%%%%%%%%%%%%%%%%%%%%%%%%%%%%%%%%%%%
%%%%%%%%%%%%%%%%%%%%%%%%%%%%%%%%%%%%%%%%%%%%%%%%%%%%%%%%%%%%%%%%%%%%%%%%%%%%%%%%%%%%%%%%%%%%%%%%%%

\section{Connected graded Calabi-Yau algebras}
\subsection{Calabi-Yau algebras, graded Calabi-Yau algebras, and AS-Gorenstein  algebras}\label{grCYASG}

Let $A$ be an associative $k$-algebra with unit. Denote by 
$A\stackrel{\mathrm{out}}{\otimes}A$ (\textit{resp.} $A\stackrel{\mathrm{inn}}{\otimes}A$) the vector space
$A\otimes A$ endowed with its outer (\textit{resp.} inner) bimodule structure. Consider the $\kk$-algebra
$A^e=A\otimes A^{op}$, so that the categories $A$-Mod-$A$, $A^e$-Mod and Mod-$A^e$ are 
naturally isomorphic. Clearly, via these isomorphisms, $A\stackrel{\mathrm{out}}{\otimes}A$ 
(\textit{resp.} $A\stackrel{\mathrm{inn}}{\otimes}A$) is identified with the left (\textit{resp.} right) $A^e$-module $A^e$
with the left (\textit{resp.} right) multiplication as action. 

The bifunctor Hom in the categories $A$-Mod-$A$, $A^e$-Mod and Mod-$A^e$ are denoted 
by $\mathrm{Hom}_{A-A}(.,.)$, $\mathrm{Hom}_{A^e-}(.,.)$, $\mathrm{Hom}_{-A^e}(.,.)$ respectively, and similarly for Ext.
The right $A^e$-module $\mathrm{Ext}_{A^{e}-}^{i}(A, A^{e})$ and the $A$-$A$-bimodule 
$\mathrm{Ext}_{A-A}^{i}(A, A\stackrel{\mathrm{out}}{\otimes}A)$ are naturally identified, as well as
$\mathrm{Ext}_{-A^{e}}^{i}(A, A^{e})$ and 
$\mathrm{Ext}_{A-A}^{i}(A, A\stackrel{\mathrm{inn}}{\otimes}A)$. If $A$ is an associative $\kk$-algebra which has a finite projective $A$-bimodule resolution by bimodules of finite type, then the left version in the 
following definition coincides with Ginzburg's definition.
\Bdf \label{cy}
We say that $A$ is a left Calabi-Yau algebra of dimension $n\pgq 1$ if 
$$\mathrm{Ext}_{A^{e}-}^{i}(A, A^{e})\cong \left\{ 
\begin{array}{ll}
A & \mbox{if }i=n \\
0 & \mbox{otherwise}
\end{array} \right. $$
in Mod-$A^e$. Replacing $A^{e}-$ by $-A^e$ and Mod-$A^e$ by $A^e$-Mod, we say that $A$ is a 
right Calabi-Yau algebra of dimension $n$.
\Edf

Since $\mathrm{Ext}_{A^{op}-A^{op}}^{i}(A^{op}, A^{op}\stackrel{\mathrm{out}}{\otimes}A^{op})\cong 
\mathrm{Ext}_{A-A}^{i}(A, A\stackrel{\mathrm{inn}}{\otimes}A)$, $A^{op}$ is left Calabi-Yau if and
only if $A$ is right Calabi-Yau. We are now interested in a graded version of Definition \ref{cy}.
\Bdf \label{grcy}
Let $A$ be a connected $\mathbb{N}$-graded $\Bbbk$-algebra. 
We say that $A$ is a left graded Calabi-Yau algebra of dimension $n\pgq 1$ if 
$$\underline{\mathrm{Ext}}_{A^{e}-}^{i}(A, A^{e})\cong \left\{ 
\begin{array}{ll}
A(\ell) & \mbox{if }i=n \\
0 & \mbox{otherwise}
\end{array} \right. $$
in grMod-$A^e$ for some $\ell \in \mathbb{Z}$ (called the parameter of $A$).
Replacing $A^{e}-$ by $-A^e$ and grMod-$A^e$ 
by $A^e$-grMod, we say that $A$ is a right graded Calabi-Yau algebra of dimension $n$.
\Edf

The following proposition is the Calabi-Yau analogue of a result by D. Stephenson and J. Zhang
concerning  graded (\textit{i.e.} Artin-Schelter) Gorenstein algebras (Proposition 3.1 
in~\cite{sz:gro}). 
\Bpo\label{ptiesgrCY}
Assume that $A$ is a left graded Calabi-Yau algebra of dimension $n$ and parameter $\ell$.
Assume that $A$ has a finite global dimension $D$. Then:
\begin{enumerate}[{\it (i)}]
\item We have $n=D$.
\item The bimodule $A$ has a graded free resolution of finite type in $A$-grMod-$A$.
\item  For any minimal projective resolution $\mathcal{P}$ of $A$ in $A^e$-grMod, 
$\underline{\mathrm{Hom}}_{A^e-}(\mathcal{P}, A^e)$ is a minimal projective resolution of
$A(\ell)$ in grMod-$A^e$.
\item $A$ is AS-Gorenstein of dimension $n$ and parameter $\ell$, and $\ell$ is nonnegative.
\item $A$ is right graded Calabi-Yau of dimension $n$ and parameter $\ell$.
\item $A$ is left and right Calabi-Yau of dimension $n$.
\end{enumerate}
\Epo
\Bdm
Let $\mathcal{P}=(P_i)_{i\pgq 0}$ be a minimal projective resolution of $A$ in $A^e$-grMod. We know that
the length of $\mathcal{P}$ equals $D$~\cite[Th\'eor\`eme 3.3]{rb:dimHoc}. In particular, $P_i=0$ if and only if $i>D$. 
For $0\ppq i \ppq D$, set $P_i = A^e\otimes E_i$, where $E_i$ is in $\Bbbk$-grMod such that for $i>0$, $E_i$ is 
a graded vector subspace of $P_{i-1}$ and the differential $d_i: P_i\rightarrow P_{i-1}$ is the natural 
$A^e$-linear extension of the inclusion $E_i \hookrightarrow P_{i-1}$.
 
Set $\underline{\mathrm{Hom}}_{A^e-}(P_i,A^e)=P_i^{\vee}$ and 
$\underline{\mathrm{Hom}}_{A^e-}(d_i,A^e)=d_i^{\vee}$. We identify $P_i^{\vee}$ with
$\underline{\mathrm{Hom}}_{\kk}(E_i,A^e)$, so that the following direct sum holds in $\Bbbk$-grMod:
\begin{equation} \label{somdir}
P_i^{\vee}= \underline{\mathrm{Hom}}_{\Bbbk}(E_i,\Bbbk) \oplus \underline{\mathrm{Hom}}_{\Bbbk}(E_i,A^e_{>0}),
\end{equation}
where $A^e_{>0}=\bigoplus_{u+v>0} A_u\otimes A_v$. The projection $p_i$ of $P_i^{\vee}$ onto the first component assigns to each $f \in P_i^{\vee}$ the linear form $\epsilon \circ f$, where 
$\epsilon:A^e\rightarrow \Bbbk$ is the 
natural projection.
Since $d_{i-1}$ is injective on $E_{i-1}$, $E_i$ is contained in $A^e_{>0}\otimes E_{i-1}$. By $A^e$-linearity of $g \in P_{i-1}^{\vee}$, $d_i^{\vee}(g)(v)=g(d_i (v))$ belongs to $A^e_{>0}$ for any $v \in E_i$. Therefore, $p_i(d_i^{\vee}(g))=0$ for any $g \in P_{i-1}^{\vee}$. Thus we have obtained
\begin{equation} \label{inc}
\im(d_i^{\vee}) \subseteq \underline{\mathrm{Hom}}_{\Bbbk}(E_i,A^e_{>0}).
\end{equation}
We have $\mathrm{Ext}_{A^{e}-}^{D}(A, A^{e})\cong P_D^{\vee}/\im(d_D^{\vee})$, and (\ref{somdir}), 
(\ref{inc}) show that this quotient maps onto $\underline{\mathrm{Hom}}_{\Bbbk}(E_D,\Bbbk)$. The latter is nonzero, 
thus $D=n$.

Note that $\underline{\mathrm{Hom}}_{\Bbbk}(E_i,A^e_{>0})$ is a sub-$A^e$-module of $P_i^{\vee}$, and the structure of
$A^e$-module inherited by $\underline{\mathrm{Hom}}_{\Bbbk}(E_i,\Bbbk)$ is trivial, \textit{i.e.} the action of $x\in A^e$ on 
$f:E_i\rightarrow \Bbbk$ is given by $f.x=f \epsilon(x)$. But now the surjective arrow 
$A\rightarrow \underline{\mathrm{Hom}}_{\Bbbk}(E_n,\Bbbk)$ of grMod-$A^e$ implies that 
$\underline{\mathrm{Hom}}_{\Bbbk}(E_n,\Bbbk)$ is generated as a vector space by the image of $1_A,$ therefore it is 1-dimensional, thus $E_n$ is 1-dimensional as well.

Let us prove that $E_i$ is finite dimensional by using a finite decreasing induction. Fix $0\ppq i\ppq D-1$ such that $E_{i+1}$ is finite dimensional and $E_i$ infinite dimensional. There exist graded vector subspaces $F$ and
 $G$ in $E_i$ such that $E_i=F\oplus G$, with $F$ finite dimensional and $E_{i+1}\subseteq A^e\otimes F$, so that
 $\im(d_{i+1})\subseteq A^e\otimes F$. Thus $(A^e\otimes F)^{\perp} \subseteq \ker(d_{i+1}^{\vee})$. But $\ker(d_{i+1}^{\vee})=\im(d_i^{\vee})$ since  $\underline{\mathrm{Ext}}_{A^e}^{i}(A,A^e)=0$ by assumption. As
 $G\neq 0$, there exists $f\in \underline{\mathrm{Hom}}_{\Bbbk}(E_i,\Bbbk)$ vanishing on $F$ and nonvanishing on $G$. Therefore $f \in (A^e\otimes F)^{\perp}\subseteq \im(d_i^{\vee})$, 
which contradicts (\ref{somdir}) and (\ref{inc}). Consequently, \textit{(i)} and \textit{(ii)} are proved. 

A consequence of \textit{(ii)} is that $P_i^{\vee}\cong E_i^{\ast}\otimes A^e$ in grMod-$A^e$, where $E_i^{\ast}$ stands for the graded $\kk$-dual. Thus $\underline{\mathrm{Hom}}_{A^e-}(\mathcal{P},A^e)$ is a projective resolution  of $A(\ell)$ in grMod-$A^e$ (using the fact that $A$ is left graded Calabi-Yau).  Another consequence of \textit{(ii)} is that any module $P_i^{\vee}$ is bounded below. In order to prove that the resolution is minimal, it is sufficient to prove that the differential of the complex $\underline{\mathrm{Hom}}_{A^e-}(\mathcal{P},A^e)\otimes _{A^e} \Bbbk$ vanishes~\cite{rb:dimHoc}. But this is straightforward from the inclusion (\ref{inc}) which can be written now as $\im(d_i^{\vee}) \subseteq E_i^{\ast}\otimes A^e_{>0}$. Hence \textit{(iii)}.

Following~\cite{rb:dimHoc}, $\mathcal{Q}=\mathcal{P}\otimes_A \Bbbk$ is a minimal projective resolution of 
$\Bbbk(\ell)$ in $A$-grMod. For the same reason, $\Bbbk \otimes_A \underline{\mathrm{Hom}}_{A^e-}(\mathcal{P},A^e)$ is a minimal 
projective resolution of $\Bbbk (\ell)$ in grMod-$A$. But 
$$\Bbbk \otimes_A \underline{\mathrm{Hom}}_{A^e-}(\mathcal{P},A^e)\cong 
\underline{\mathrm{Hom}}_{A-}(\mathcal{Q},A)$$
in grMod-$A$. In particular, $\underline{\mathrm{Hom}}_{A-}(\mathcal{Q},A)$ is a resolution of $\Bbbk (\ell)$, so that 
$$\underline{\mathrm{Ext}}_{A-}^{i}(\Bbbk, A)\cong \left\{ 
\begin{array}{ll}
\Bbbk(\ell) & \mbox{if }i=n \\
0 & \mbox{otherwise}
\end{array} \right. $$
in grMod-$A$. Thus $A$ is AS-Gorenstein of dimension $n$, and the Gorenstein parameter is $\ell$. It is known that $\ell \pgq 0$~\cite{sz:gro}. Hence \textit{(iv)}.

To prove \textit{(v)}, we apply the functor $\underline{\mathrm{Hom}}_{-A^e}(\cdot ,A^e)$ to the 
projective resolution $\underline{\mathrm{Hom}}_{A^e-}(\mathcal{P},A^e)$ of $A(\ell)$. Since the 
$E_i$'s are finite dimensional, we recover the resolution $\mathcal{P}$ of $A$. Thus 
$$\underline{\mathrm{Ext}}_{-A^{e}}^{i}(A(\ell), A^{e})\cong \left\{ 
\begin{array}{ll}
A & \mbox{if }i=n \\
0 & \mbox{otherwise}
\end{array} \right. $$
Shifting by $\ell$, we get \textit{(v)}. 

A consequence of \textit{(ii)} is that 
$$\mathrm{Ext}_{A^{e}-}^{i}(A, A^{e}) \cong \underline{\mathrm{Ext}}_{A^{e}-}^{i}(A, A^{e})$$
in Mod-$A^e$, and the same on the right. Hence \textit{(vi)}.
\Edm

\br\label{CYgrCY}
Let $A$ be a connected $\mathbb{N}$-graded $\Bbbk$-algebra. Assume that \textit{(ii)} holds. If $A$ is left 
Calabi-Yau of dimension $n$, then $A$ is left graded Calabi-Yau of dimension $n$.
\er

%%%%%%%%%%%%%%%%%%%%%%%%%%%%%%%%%%%%%%%%%%%%%%%%%%%%%%%%%%%%%%%%%%%%%%%%%%%%%%%%%%%%%%%%%%%%%%%%%%%%%%%%%%%%%%%%%%%%%%%%%%%%%%%%%%%%%%%%%%%%%%%%%%%%%%%%%%%%%%%%%%%

\subsection{A criterion for certain graded algebras to be Calabi-Yau}\label{criterion}

This criterion (Proposition below) will be useful to us for Examples \ref{yangmills}, \ref{ASregular} and \ref{antisymmetriser}. 
Let $A$ be a connected $\mathbb{N}$-graded $\Bbbk$-algebra. Let us assume that $A$ is $N$-Koszul 
for a certain $N\pgq 2$, and that $A$ is AS-Gorenstein of global dimension $n$. It is 
known~\cite[Section 6]{rbnm:kogo} that $\mathrm{Ext}_{A^{e}-}^{i}(A, A^{e})=0$ if $i\neq n$, and that
$$\mathrm{Ext}_{A^{e}-}^{n}(A, A^{e})\cong_{\varepsilon^{n+1}\phi}A(\zeta _N(n))$$ in 
the category $A$-grMod-$A$ of graded bimodules, where $\zeta _N$ 
stands for the jump function of the $N$-Koszul algebra $A$. Here $\varepsilon $ and $\phi $ are 
certain automorphisms of the graded algebra $A$. We refer to~\cite{rbnm:kogo} for more details.
\Bpo \label{auto}
With the notations and assumptions above, the algebra $A$ is Calabi-Yau 
if and only if $\varepsilon^{n+1}\phi=\mathrm{id}_A$.
\Epo
\Bdm
It is clear from Ginzburg's definition that, if $\varepsilon^{n+1}\phi=\mathrm{id}_A$, 
$A$ is Calabi-Yau of dimension $n$. Conversely, assume that $A$ is Calabi-Yau (necessarily 
of dimension $n$). Remark \ref{CYgrCY} and Proposition \ref{ptiesgrCY}  show
that $A$ is graded Calabi-Yau of dimension $n$, and that the Calabi-Yau parameter coincides with 
the Gorenstein parameter $\zeta _N(n)$. Thus there is an isomorphism 
$f: A\rightarrow _{\varepsilon^{n+1}\phi}A$ in $A$-grMod-$A$, so that $f(1)$ is a nonvanishing 
element in $\Bbbk$. For any $a$, $b$ in $A$, we have $f(ab)=\varepsilon^{n+1}\phi(a)f(1)b$, 
implying that $f(a)=\varepsilon^{n+1}\phi(a)f(1)=f(1)a$, and $\varepsilon^{n+1}\phi(a)=a$ follows.
\Edm

%%%%%%%%%%%%%%%%%%%%%%%%%%%%%%%%%%%%%%%%%%%%%%%%%%%%%%%%%%%%%%%%%%%%%%%%%%%%%%%%%%%%%%%%%%%%%%%%%%%%%%%%%%%%%%%%%%%%%%%%%%%%%%%%%%%%%%%%%%%%%%%%%%%%%%%%%%%%%%%%%%%%%%%%%%%%%%%%%%%%%%%%%%%%%%%%%%%%%%%%%%%%%%%%%%%%%%%%%%%%%%%%%%%%%%%%%%%%%%%%%%%%%%%%%%%%%%%%%%%%%%%%%%%%%%%%%%%%%%%%%%%%%%%%%%%%%%%%%%%%%%%%%%%%%%%%%%%%%%%%%%%%%%%%%%%%%%%%%%%%%%%%%%%%%%%%%%%%%%%%%%%%%%%%%%%%%%%%%%%%%%%%%%%%%%%%
%%%%%%%%%%%%%%%%%%%%%%%%%%%%%%%%%%%%%%%%%%%%%%%%%%%%%%%%%%%%%%%%%%%%%%%%%%%%%%%%%%%%%%%%%%%%%%%%%%

\section{Examples}\label{examples}

Throughout Examples \ref{yangmills}, \ref{ASregular} and \ref{antisymmetriser}, the characteristic of $\kk$ is assumed to be zero.

\bex\label{yangmills} \textbf{Yang-Mills algebra.}

The Yang-Mills algebra $A$ with $s+1$ generators is graded, 3-Koszul, AS-Gorenstein of global 
dimension 3~\cite{cdv:yang}. Following~\cite{cdv:yang2}, the identity
$\varepsilon^{n+1}\phi=\mathrm{id}_A$ holds, thus $A$ is Calabi-Yau by Proposition \ref{auto} 
(equivalently, it would be possible to extract from~\cite{cdv:yang2} a self-dual resolution
of $A$ by bimodules, and then $A$ is Calabi-Yau by Theorem 4.2 in~\cite{bock:gcy}).

Let us show directly that $A$ derives from a potential (see also~\cite{vg:project}). 
We use the material contained in~\cite{rbdv:iym}. The generators of the 
Yang-Mills algebra $A$ are denoted by $\nabla_0, \nabla_1, \ldots , \nabla_s$, where $s\pgq 1$. 
Let $(g^{\alpha \beta })_{0\ppq \alpha ,\beta \ppq s}$ be a nondegenerate symmetric 
matrix with entries in $\Bbbk$. The relations of $A$ are denoted by $W^0, W^1, \ldots , W^s$, 
where $W^{\rho }=\sum_{\lambda ,\mu ,\nu }
W^{\rho \lambda \mu \nu } \nabla_{\lambda }\nabla_{\mu }\nabla_{\nu }$, and
$$W^{\rho \lambda \mu \nu }=g^{\rho \lambda }g^{\mu \nu } + g^{\rho \nu }g^{\lambda \mu} 
- 2 g^{\rho \mu }g^{\lambda \nu }.$$
Using the concepts introduced in~\cite{mdv:multi,dubois-violette}, $A$ is defined from the 
following "volume form"
$$W_4=\sum_{\rho ,\lambda ,\mu ,\nu }
W^{\rho \lambda \mu \nu } \nabla_{\rho }\nabla_{\lambda }\nabla_{\mu }\nabla_{\nu }.$$
Actually, this volume form $W_4$ is a potential because of the 
relation $W^{\rho \lambda \mu \nu }=W^{\lambda \mu \nu \rho }$. Moreover we have
$\partial_{\rho }W_4= 4 W^{\rho }$ for any $\rho $. In other words, $A=A(Q,W_4)$, where $Q$ is the quiver having one vertex and $s+1$ loops.

Following~\cite{rbdv:iym}, a PBW deformation $U$ of $A$ satisfying (PBW2') has the following 
relations
$$W^{\rho } + \sum_{\mu ,\nu } c^{\mu \nu \rho }\nabla_{\mu }\nabla_{\nu } + 
\sum_{\lambda } s^{\lambda \rho }\nabla_{\lambda } + s^{\rho }=0, \ \ 0\ppq \rho \ppq s ,$$
where the coefficients are elements of $\Bbbk$ such that
$c^{\mu \nu \rho }=c^{\rho \mu \nu }, \ s^{\lambda \rho }=s^{\rho \lambda }$. Then it is easy to 
verify that $U$ derives from a potential, as stated in  Theorem \ref{NSCpbw2'}. In fact, setting
$$W_3=\frac{1}{3} \sum_{\mu ,\nu , \rho} c^{\mu \nu \rho }\nabla_{\mu }\nabla_{\nu }\nabla_{\rho },
\ W_2=\frac{1}{2} \sum_{\lambda \rho} s^{\lambda \rho }\nabla_{\lambda }\nabla_{\rho }, 
\ W_1=\sum_{\rho }s^{\rho }\nabla_{\rho },$$
we have $U=A(Q, W_4 + W_3 + W_2 + W_1)$. 

All the PBW deformations of $A$ are determined in~\cite{rbdv:iym}. From these 
computations, it is immediate to see that there exist PBW deformations of $A$ for which (PBW2') 
does not hold.

\br
Let us show how we can deduce the formula $\varepsilon^{n+1}\phi=\mathrm{id}_A$
from the definitions of $\varepsilon $ and $\phi $ contained in~\cite{rbnm:kogo}. Since $n=3$, 
$\varepsilon^{n+1}=\mathrm{id}_A$ is clear. 
The automorphism $\phi $ is uniquely determined by its component $\phi _1$, and $\phi _1$
is defined as the transposed linear map of $\nu _1$ where $\nu $ is the 
Nakayama automorphism of the Yoneda algebra $E(A)$ of $A$. Our aim is now 
to prove that $\nu _1=\mathrm{id}_{E(A)_1}$. We have to compute the Frobenius 
pairing $(.,.)$ of $E(A)$. Actually, it suffices to compute $(x,y)$ and $(y,x)$ when $x\in E(A)_1=V^{\ast}$ 
and $y\in E(A)_2=R^{\ast}$. As usual, $V$ denotes the generator space of $A$
and $R$ denotes the relation space. We have (formula (6.3) in~\cite{rbnm:kogo})
$$(x,y)=\langle x\otimes g,W_4\rangle , \ \ (y,x)=\langle g\otimes x,W_4\rangle ,$$
where $g$ is any linear form on $V^{\otimes 3}$ extending $y:R\rightarrow \Bbbk$, and 
consequently $x\otimes g$ and $g\otimes x$ are linear forms on $V^{\otimes 4}$. Then $\nu _1$
is defined by $(x,y)=(y,\nu _1(x))$.

Since $W_4=\sum_{\alpha }\nabla_{\alpha }W^{\alpha }$, we have 
$(\nabla_{\lambda }^{\ast},W^{\alpha \ast}) = \delta _{\lambda \alpha }$ (Kronecker symbol).
Similarly, $W_4=\sum_{\alpha }W^{\alpha }\nabla_{\alpha }$ implies
$(W^{\alpha \ast},\nabla_{\lambda }^{\ast}) = \delta _{\lambda \alpha }$. Thus 
$\nu _1=\mathrm{id}_{E(A)_1}$.
\er
\eex

\bex\label{ASregular} \textbf{Artin-Schelter regular algebras of global dimension $3$.}

In this example, $k$ is the field of complex numbers. Applying Proposition \ref{auto} in 
the quadratic case~\cite[Corollary 9.3]{vdb:exist} and in the cubic 
case~\cite[Proposition 6.5]{rbnm:kogo}, we get the following.
\Bpo\label{typeofASregular}
Let $A$ be an AS-regular algebra (with polynomial growth) of global dimension 3. Then 
$A$ is Calabi-Yau if and only if $A$ is of type {\rm A} in the classification of Artin and 
Schelter~\cite{as:regular}.
\Epo
It is easy to show directly that $A$ of type A derives from a potential. Let us limit ourselves 
to the cubic case, and let us examine the PBW deformations in this case. If $A$ is cubic of type A, 
then $A=A(Q, W_4)$ where the quiver $Q$ has one vertex and two loops $x$ and $y$, and the
potential is
$$W_4=y^4 + a(x^2y^2+xy^2x+y^2x^2+yx^2y) + b(xyxy+yxyx) + x^4.$$
We have $\partial_x W_4 = 4f$ and $\partial_y W_4 = 4g$, where
$$f=ay^2x+byxy+axy^2+x^3,\ \ g=y^3+ayx^2+bxyx+ax^2y$$
are the relations of $A$. Here $a$ and $b$ are the usual generic complex parameters.

In~\cite{fv:pbwdef}, Fl\o ystad and Vatne have determined all the PBW deformations of $A$ 
(actually for any type, not necessarily  type A). From their computations, 
it results that any PBW deformation $U$
of $A$ (of type A) satisfies condition (PBW2'). Moreover $U$ has
the following relations 
$$f + a_{11}x^2 + b_{11}xy + b_{11}yx + a_{14}y^2 + a_{21}x + a_{22}y + a_3 = 0,$$
$$g + b_{11}x^2 + a_{14}xy + a_{14}yx + b_{14}y^2 + a_{22}x + b_{22}y + b_3 = 0,$$
where the coefficients are arbitrary complex numbers. Then it is easy to 
verify that $U$ derives from a potential, as stated in Theorem \ref{NSCpbw2'}. In fact, setting
\begin{align*}
&W_3=\frac{1}{3} \left(a_{11}x^3 + b_{11}(x^2y+xyx+yx^2) + a_{14}(y^2x+yxy+xy^2) + 
b_{14}y^3\right),\\
&W_2=\frac{1}{2} \left(a_{21}x^2 + a_{22}(xy+yx) + b_{22}y^2 \right),\\
&W_1=a_3x + b_3y,
\end{align*}
we have $U=A(Q, W_4 + W_3 + W_2 + W_1)$. 
\eex

\bex\label{antisymmetriser} \textbf{Antisymmetriser algebras.}

Let us introduce some general notations. If $p$ is a cycle of a quiver $Q$, denote by 
$\bar{p}$ the class of $p$ up to cyclic permutation of the arrows. Denote by $\mathcal{C}$ the 
set of cycles of $Q$ and by $\overline{\mathcal{C}}$ the set of classes of cycles in $\C$ modulo cyclic permutations. The natural map from
$\mathcal{C}$ into the set $\mathrm{Pot}(Q)$ of potentials of $Q$ induces a linear isomorphism
from $\Bbbk \overline{\mathcal{C}}$ to $\mathrm{Pot}(Q)$. We shall often identify  $\kk \ovl{\C}$ and $\pot{Q}$  via this isomorphism, and $\overline{\mathcal{C}}$ will be considered 
as a basis of the vector space $\mathrm{Pot}(Q)$.

Now assume that $\Bbbk$ has characteristic zero, $Q$ has one vertex and $n\pgq 2$ loops $x_1,\ldots, x_n$. 
Any path is a cycle.
Introduce $W_n=\mathrm{Ant}(x_1,\ldots, x_n)$ where $\mathrm{Ant}$ stands for the antisymmetriser
of the variables on which this operator is acting, \textit{i.e.}
$$W_n= \sum_{\sigma \in S_n} \mathrm{sgn} (\sigma )\, x_{\sigma (1)} \ldots x_{\sigma (n)}.$$
Here $S_n$ denotes the group of permutations of $1,\ldots ,n$.
\Blm \label{lemant1}
If $n$ is even, $\overline{W}_n=0$. If $n$ is odd, then for any $i=1,\ldots ,n$ we have
the formulas
$$\overline{W}_n=n \, \sum_{\sigma (1)=i} \mathrm{sgn} (\sigma ) \ 
\overline{x_{\sigma (1)} \ldots x_{\sigma (n)}},$$
$$\partial _{x_i}\overline{W}_n=n \, (-1)^{i+1} \mathrm{Ant}(x_1,\ldots, \widehat{x_i}, \ldots, x_n),$$
in which $\widehat{x_i}$ means that the variable $x_i$ is removed.
\Elm
\Bdm
Fix $i$ and write 
\begin{align*}
& W_n  = \sum_{\sigma (1)=i} \mathrm{sgn} (\sigma )  (x_{\sigma (1)} \ldots x_{\sigma (n)} + 
    (-1)^{n+1} x_{\sigma (2)} \ldots x_{\sigma (n)} x_{\sigma (1)} \\
& + (-1)^{2(n+1)} x_{\sigma (3)} \ldots x_{\sigma (n)} x_{\sigma (1)} x_{\sigma (2)} + \cdots 
    + (-1)^{(n-1)(n+1)} x_{\sigma (n)} x_{\sigma (1)} \ldots x_{\sigma (n-1)})
\end{align*}
hence
$$\overline{W}_n=\sum_{\sigma (1)=i} \mathrm{sgn} (\sigma ) (1+(-1)^{n+1}+(-1)^{2(n+1)}+
\cdots (-1)^{(n-1)(n+1)}) \overline{x_{\sigma (1)} \ldots x_{\sigma (n)}}.$$
Clearly, $\overline{W}_n=0$ if $n$ is even, and we have the first formula if $n$ is odd.
From this first formula, we deduce 
$$\partial _{x_i}\overline{W}_n=n \sum_{\sigma (1)=i} \mathrm{sgn} (\sigma )  
\ x_{\sigma (2)} \ldots x_{\sigma (n)},$$
and the second formula follows, since the sign of the permutation $1\mapsto i, 
2\mapsto 1, \ldots , i\mapsto i-1$ is $(-1)^{i+1}$.
\Edm 
\\ 

For the remainder of this example, we assume $n\pgq 3$. Denote by $A$ the $(N=n-1)$-homogeneous
algebra whose generators are $x_1,\ldots, x_n$ and relations are 
$\mathrm{Ant}(x_1,\ldots, \widehat{x_i}, \ldots, x_n)$ for $i=1,\ldots ,n$, and call it the
\textit{antisymmetriser algebra} (it was introduced in~\cite{rb:nonquad} for any $2\ppq N\ppq n$). 
Lemma \ref{lemant1} shows that $A=A(Q,W_n)$ if $n$ is odd. When $n$ is even, we can improve the statement $\overline{W}_n=0$.
\Blm \label{lemant2}
If $n$ is even, there exists no potential $W$ such that $A=A(Q,W)$.
\Elm
\Bdm
Assume that $W$ exists. Write $W=\sum_{\overline{p}\in \overline{\mathcal{C}}_n} 
\lambda _{\overline{p}} \,\overline{p}$, where $\mathcal{C}_n$ denotes the set of cycles of length $n$. Let $W'$ be the part
of $W$ for which the summation is performed only over the cycles $p$ containing (at least) 
twice the same loop. The fact that $\partial_{x_i} W$ belongs to 
$\Bbbk.\mathrm{Ant}(x_1,\ldots, \widehat{x_i}, \ldots, x_n)$ implies $\partial_{x_i} W'=0$, 
and this holds for $i=1,\ldots ,n$. 

Now it is easy to see in general that $$\cap_{a\in Q_1}\ker (\partial_a)=\kk Q_0.$$ Here $W'\in \cap_{i=0}^n\ker\partial_{x_i}$ and $W'$ is homogeneous of length $n\pgq 3$, so $W'=0.$
Therefore we can write
$$W= \sum_{\overline{\sigma} \in \overline{S}_n} \lambda _{\overline{\sigma}} \ 
\overline{x_{\sigma (1)} \ldots x_{\sigma (n)}},$$
with obvious notations for $\overline{\sigma}$ and $\overline{S}_n$. Fixing $i$, we have a 
unique decomposition
$$W= \sum_{\sigma (1)=i} \lambda _{\sigma}^{(i)} \ 
\overline{x_{\sigma (1)} \ldots x_{\sigma (n)}},$$
since in each class $\overline{\sigma}$ there is a unique permutation mapping 1 to $i$. The superscript $(i)$ in $\lambda_\sigma^{(i)}$ is there to remind us that $\sigma$ depends on $i$ (since $\sigma(1)=i$). Therefore, we obtain
$$\partial _{x_i} W=\sum_{\sigma (1)=i} \lambda _{\sigma}^{(i)} 
\ x_{\sigma (2)} \ldots x_{\sigma (n)}.$$
The fact that $\partial_{x_i} W$ belongs to 
$\Bbbk.\mathrm{Ant}(x_1,\ldots, \widehat{x_i}, \ldots, x_n)$ implies that 
$$\partial _{x_i} W = \lambda _{\tau _i}^{(i)}\ 
\mathrm{Ant}(x_{i+1},\ldots, x_n, x_1, \ldots, x_{i-1}),$$
where $\tau _i$ is the permutation 
$$1\mapsto i, 2 \mapsto i+1, \ldots , n-i+1 \mapsto n,\, n-i+2 \mapsto 1, \ldots, n\mapsto i-1.$$
Decompose $\sigma =\sigma '\tau _i$, $\sigma '$ leaving $i$ fixed. The comparison of the two 
previous decompositions of $\partial _{x_i} W$ provides $\lambda _{\sigma}^{(i)}=
\lambda _{\tau _i}^{(i)}\, \mathrm{sgn} (\sigma ')$, and we get
\begin{equation} \label{crucial}
W= \lambda _{\tau _i}^{(i)} \sum_{\sigma '(i)=i} \mathrm{sgn} (\sigma ') \ 
\overline{x_{\sigma '(i)} \ldots x_{\sigma '(n)} x_{\sigma '(1)}\ldots x_{\sigma '(i-1)}}.
\end{equation}
For $\sigma '=\mathrm{id}$ in (\ref{crucial}), we see that $\lambda _{\tau _i}^{(i)}$ coincides with
the coefficient $\lambda $ of $\overline{x_1 \ldots x_n}$ in the linear decomposition of $W$ 
in the basis $\overline{\mathcal{C}}_n$. The case $i=1$ in (\ref{crucial}) gives
$$W=\lambda [\overline{x_1 \ldots x_n} - \overline{x_1 x_3 x_2 x_4 \ldots x_n}+ \cdots ],$$
while the case $i=2$ gives
\begin{align*}
W&= \lambda [\overline{x_2x_3 \ldots x_nx_1} +(-1)^{n-2}\, \overline{x_2 x_4 \ldots x_nx_1x_3}+ \cdots ]\\
 &= \lambda [\overline{x_1 \ldots x_n} + \overline{x_1 x_3 x_2 x_4 \ldots x_n}+ \cdots ].
\end{align*}
The comparison of the two so-obtained linear decompositions in the basis $\overline{\mathcal{C}}_n$ 
implies $\lambda =0$, hence $W=0$ and we have a contradiction.
\Edm

For any $n\pgq 3$, it is known that the antisymmetriser algebra $A$ is $(n-1)$-Koszul of global
dimension 3~\cite[Theorem 3.13]{rb:nonquad}, AS-Gorenstein~\cite[Corollary 5.10]{rbnm:kogo}. Moreover the automorphism
$\varepsilon^{n+1}\phi$ is $\mathrm{id}_A$ or $-\mathrm{id}_A$ according to  whether $n$ is 
odd or even~\cite[end of the paper]{rbnm:kogo}. Thus Proposition \ref{auto} shows that $A$ is Calabi-Yau if and only 
if $n$ is odd.

Fl\o ystad and Vatne have determined the PBW deformations of any antisymmetiser algebra $A$ such 
that $N\ppq n-2$ (Theorem 4.1 and Theorem 4.2 in~\cite{fv:pbwdef}). The case $N=n-1$ of interest for
us is not treated in~\cite{fv:pbwdef}. However  Theorems \ref{nonhomogisPBW} and \ref{NSCpbw2'}  provide the following 
partial answer to the determination of PBW deformations of $A$ when $n\pgq 3$ is odd and $N=n-1$:
\begin{enumerate}[(a)] 
\item \label{antisymm1} For any $W_j\in \mathrm{Pot}(Q)_j$, $1\ppq j \ppq n-1$, the algebra $A(Q, \sum_{1\ppq j\ppq n} W_j)$ is 
a PBW deformation of $A=A(Q,W_n)$.
\item \label{antisymm2} Any PBW deformation in (\ref{antisymm1}) satisfies (PBW2').
\item \label{antisymm3}  Any PBW deformation of $A$ satisfying (PBW2') is as in (\ref{antisymm1}).
\end{enumerate}

This leads to the following questions:
\begin{enumerate}[{Question} 1.]
\item Determine $\dim (\mathrm{Pot}(Q)_j)$ for $1\ppq j \ppq n-1$.
\item Does (PBW2') hold for any PBW deformation of $A$?
\end{enumerate}

We examine now Questions 1 and 2 in the simplest case $n=3$. In this case, $A=k[x,y,z]$. The 
answer to Question 1 is easy since $\overline{\mathcal{C}}_1=\{x,y,z\}$ and 
$\overline{\mathcal{C}}_2=\{x^2,y^2,z^2,xy,yz,zx\}$. Let us examine the \textit{linear}
PBW deformations of $A$, \textit{i.e.} such that $\varphi _0=0$. It is well known that the 
linear PBW deformations of the polynomial algebra $A$ are exactly 
the Lie algebras having $V=\Bbbk x\oplus \Bbbk y\oplus \Bbbk z$ as underlying vector space (see \textit{e.g.} Example 3.7 
in~\cite{rbvg:hisym}). More precisely, let $\varphi _1:R\rightarrow V$ be linear, and denote by 
$[.,.] :V\times V \rightarrow V$ the antisymmetric bilinear map extending naturally $\varphi _1$ 
($R$ is spanned by $yz-zy, zx-xz, xy-yx$).
Then $\varphi _1$ is a PBW deformation of $A$ if and only if the relation
$$\varphi _1(\varphi _1\otimes \mathrm{id}_V - \mathrm{id}_V \otimes \varphi _1) 
(W_3)=0$$
holds, which in turn is equivalent to the Jacobi identity of the bracket $[.,.]$. 

Let us give the matrix $\alpha =(\alpha _{ij})_{1\ppq i,j\ppq 3}$ of coefficients of 
$\varphi _1$, that is 
\begin{align*}
& \varphi _1 (yz-zy) = \alpha _{11}x + \alpha _{12}y + \alpha _{13}z,\\
& \varphi _1 (zx-xz) = \alpha _{21}x + \alpha _{22}y + \alpha _{23}z,\\
& \varphi _1 (xy-yx) = \alpha _{31}x + \alpha _{32}y + \alpha _{33}z.
\end{align*}
Then we obtain
\begin{equation}\begin{array}{rl} \label{pbw2}
(\varphi _1\otimes \mathrm{id}_V - \mathrm{id}_V \otimes \varphi _1)(W_3)=&
 (\alpha _{32}-\alpha _{23})\,(yz-zy)+(\alpha_{13}-\alpha_{31})\,(zx-xz)\\&+(\alpha_{21}-\alpha_{12})\,(xy-yx),
\end{array}\end{equation}
and the Jacobi identity is equivalent to the algebraic system formed by the equation
$$\alpha _{12}\alpha _{31} - \alpha _{13}\alpha _{21} + 
\alpha _{11} (\alpha _{23}-\alpha _{32})=0,$$
and the two other equations deduced by cyclic permutations of indices. The following is immediate from
(\ref{pbw2}).
\Blm \label{pbw2'}
Condition (PBW2') is equivalent to saying that the matrix $\alpha $ is symmetric. In this case, the potential is given by $W_2=-\frac{1}{2}\alpha_{11}x^2-\frac{1}{2}\alpha_{22}y^2-\frac{1}{2}\alpha_{33}z^2-\alpha_{12}xy-\alpha_{23}yz-\alpha_{31}zx.$
\Elm

Thus the answer to Question 2 is no, since it is easy to find $\alpha $ non symmetric 
and satisfying Jacobi. Actually, it is possible to be more precise. Let us fix an arbitrary
element $r=a(yz-zy)+b(zx-xz)+c(xy-yx)$ in $R$. Then, there exists a unique \textit{antisymmetric} 
$\alpha $ (denoted by $\alpha _0$) satisfying Jacobi and 
$(\varphi _1\otimes \mathrm{id}_V - \mathrm{id}_V \otimes \varphi _1)(W_3)=r$. Furthermore, 
the set of these $\alpha $'s (not necessarily antisymmetric) is equal to
$$\set{\beta + \alpha _0\,\mid\, \beta \mbox{ symmetric and } \beta\mx{a\\b\\c}=0 }$$
which is a linear affine space of dimension 6 if $a=b=c=0$, of dimension
3 otherwise.

It would be interesting to have an analogue of the above discussion for higher $n$.

\eex

We now consider quiver algebras with several vertices, described in \cite{bock:gcy}.

\bex Set $N=\ell k-1$ with $\ell\pgq2$ and $k\pgq3$, and let $A=A(Q,W_{N+1})$ be the Calabi-Yau algebra of dimension $3$ defined by the quiver $Q:$
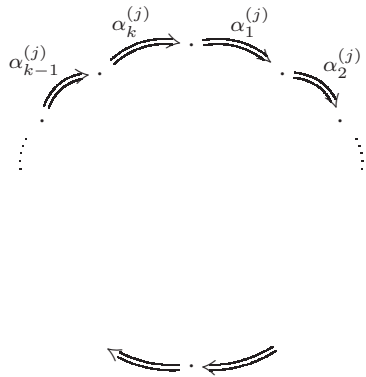
$$\xymatrix@=.01cm{
&&&&&&&&&&&\cdot\ar@/^.5pc/@{=>}[rrrrrd]^{\alpha^{(j)}_1}\\
&&&&&&\cdot\ar@/^.5pc/@{=>}[rrrrru]^{\alpha^{(j)}_{k}}&&&&&&&&&&\cdot\ar@/^.5pc/@{=>}[rrrdd]^{\alpha^{(j)}_2}\\\\
&&&\cdot\ar@/^.5pc/@{=>}[rrruu]^{\alpha^{(j)}_{k-1}}\ar@{.}@/_.3pc/[lddd]&&&&&&&&&&&&&&&&\cdot\ar@{.}@/^.3pc/[rddd]\\\\
&&&&&&&&&&&&&&&&&&&&&&&\\
&&&&&&&&&&&&&&&&&&&&&&&\\
\\
\\\\\\\\\\\\\\\\
&&&&&&&&&&&&&&&&\\
&&&&&&&&&&&\cdot\ar@/_.3pc/@{.}@{<=}[rrrrru]\ar@/^.3pc/@{.}@{=>}[lllllu]
}$$ containing $k$ vertices $1,\ldots,k$ and $n_i$ arrows $\alpha^{(1)}_i,\ldots,\alpha^{(n_i)}_i$  from the vertex $i$ to the vertex ${i+1}$ with $n_i\pgq 2$ for all $i$ with $1\ppq i\ppq k$, and by the potential 
\begin{align*}
W_{N+1}=&\sum_{i\in Q_0} \alpha_i^{(1)}\alpha_{i-1}^{(1)}\alpha_{i-2}^{(2)}\ldots\alpha_{i-N}^{(2)}
\\&+\sum_{\tiny
\begin{array}{l}
c\in Q_1\\c\not\in\set{\alpha_i^{(1)},\alpha_i^{(2)}}
\end{array}
}c\alpha_{t(c)-1}^{(1)}\alpha_{t(c)-2}^{(2)}\ldots\alpha_{t(c)-k+1}^{(2)}(c\alpha_{t(c)-1}^{(2)}\alpha_{t(c)-2}^{(2)}\ldots\alpha_{t(c)-k+1}^{(2)})^{\ell-1}.
\end{align*} We assume that the characteristic of $\kk$ does not divide $N!.$

Now let $A'$ be a PBW deformation of $A.$ Then $A'$ satisfies the following condition: $$(\star)\  \forall e\in Q_0,\ \sum_{b\in Q_1e}\varphi_{N-1}(\partial_bW_{N+1})b-\sum_{a\in eQ_1}a\varphi_{N-1}(\partial_aW_{N+1})\in\mathrm{span}_\kk\set{\partial_cW_{N+1}\,\mid\,c\in Q_1}.$$ 

However, $\varphi_{N-1}(\partial_\alpha W_{N+1})$ has length $N-1=\ell k-2,$ starts at $t(\alpha)$ and ends at $s(\alpha),$ so that $\ell k-2$ must be of the form $\lambda k-1$, which is impossible since $k\pgq 3.$ Therefore $\varphi_{N-1}(\partial_\alpha W_{N+1})=0$ for all $\alpha.$ In particular, in the relation $(\star)$, the left-hand side is equal to $0,$ and therefore $A'$ is defined by a potential (since (PBW2') holds).

Note that, in the same way as above, we can check that $\varphi_{j-1}(\partial_\alpha W_{N+1})=0$ whenever $j\not\equiv 0 \pmod{k}.$

Let us now give a more specific example: assume that $n_i=2$ for every $i\in Q_0$ (there are exactly two arrows between consecutive vertices), so that $W_{N+1}=\sum_{i\in Q_0} \alpha_i^{(1)}\alpha_{i-1}^{(1)}\alpha_{i-2}^{(2)}\ldots\alpha_{i-N}^{(2)}$. Let $A'$ be defined by $\varphi_{\lambda k-1}(\partial_{\alpha_i^{(1)}}W_{N+1})=\alpha_{i+\lambda k-1}^{(1)}\ldots \alpha_{i+2}^{(1)}\alpha_{i+1}^{(1)}$ and $\varphi_{\lambda k-1}(\partial_{\alpha_i^{(2)}}W_{N+1})=\alpha_{i+\lambda k-1}^{(2)}\ldots \alpha_{i+2}^{(2)}\alpha_{i+1}^{(2)}$. Then the potential defining $A'$ is $W=\sum_{p=1}^{N+1}W_p$ with $W_p=0$ if $p\not\equiv 0 \pmod{k},$ and $W_{\lambda k}=\sum_{i\in Q_0}\left(\alpha_{i+\lambda k-1}^{(1)}\ldots \alpha_{i+1}^{(1)}\alpha_{i}^{(1)}+\alpha_{i+\lambda k-1}^{(2)}\ldots \alpha_{i+1}^{(2)}\alpha_{i}^{(2)}\right).$
  
\eex

\bex Let $N$ be an integer with $N\pgq 3$,  and let $A=A(Q,W_{N+1})$ be the Calabi-Yau algebra defined by the quiver $Q:$
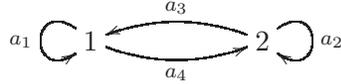
$$\xymatrix@=.05cm{1\ar@(ul,dl)[]_{a_1}\ar@/_.6pc/[rrrrrrrr]_{a_4}&&&&&&&&2\ar@/_.6pc/[llllllll]_{a_3}\ar@(ur,dr)[]^{a_2}}$$ and by the potential $W_{N+1}=a_3a_4a_1^{N-1}+a_4a_3a_2^{N-1}.$

Assume that the characteristic of $\kk$ does not divide $(N+1)!$, and let $A'$ be a PBW deformation of $A.$ 
 Then $A'$ satisfies the following conditions: $$(*_1) \varphi_{N-1}(\partial_{a_1}W_{N+1})a_1+\varphi_{N-1}(\partial_{a_4}W_{N+1})a_4-a_1\varphi_{N-1}(\partial_{a_1}W_{N+1})-a_3\varphi_{N-1}(\partial_{a_3}W_{N+1})=\mu \partial_{a_1}W_{N+1}$$ $$(*_2) \varphi_{N-1}(\partial_{a_2}W_{N+1})a_2+\varphi_{N-1}(\partial_{a_3}W_{N+1})a_3-a_2\varphi_{N-1}(\partial_{a_2}W_{N+1})-a_4\varphi_{N-1}(\partial_{a_4}W_{N+1})=\nu \partial_{a_2}W_{N+1}$$ with $\mu,\nu\in \kk $ (recall that $\varphi_{j-1}(\partial_aW_{N+1})$ is a linear combination of paths from $t(a)$ to $s(a)$).

The only terms occuring in the $\varphi_{N-1}(\partial_{a_i}W_{N+1})$ which give terms in $\partial_{a_1}W_{N+1}=\sum_{u+v=N-2}a_1^ua_3a_4a_1^v$ and in $\partial_{a_2}W_{N+1}=\sum_{u+v=N-2}a_2^ua_4a_3a_2^v$ are given as follows: 
$$
\begin{array}{l}
\varphi_{N-1}(\partial_{a_1}W_{N+1})=\sum_{u+v=N-3}\lambda_{1,u}a_1^ua_3a_4a_1^v+X_1 \\
\varphi_{N-1}(\partial_{a_2}W_{N+1})=\sum_{u+v=N-3}\lambda_{2,u}a_2^ua_4a_3a_2^v+X_2 \\
\varphi_{N-1}(\partial_{a_3}W_{N+1})=\sigma_3a_2^{N-2}a_4+\lambda_3a_4a_1^{N-2}+X_3 \\
\varphi_{N-1}(\partial_{a_4}W_{N+1})=\sigma_4a_1^{N-2}a_3+\lambda_4a_3a_2^{N-2}+X_4 
\end{array}
$$ with $\lambda_{1,u}, \lambda_{2,u}, 
\lambda_3,\lambda_4,\sigma_3$ and $\sigma_4$ in $\kk ,$ and where $X_1,X_2,X_3,X_4$ are linear combinations of paths of length $N-1$ starting and ending at appropriate vertices.

The relation $(*_1)$ gives $$\lambda_4-\lambda_{1,0}=\lambda_{1,0}-\lambda_{1,1}=\cdots=\lambda_{1,p}-\lambda_{1,p+1}=\cdots=\lambda_{1,N-3}-\sigma_3=\mu$$ and $\lambda_3=\sigma_4$. The relation $(*_2)$ gives  $$\lambda_3-\lambda_{2,0}=\lambda_{2,0}-\lambda_{2,1}=\cdots=\lambda_{2,p}-\lambda_{2,p+1}=\cdots=\lambda_{2,N-3}-\sigma_4=\nu$$ and $\lambda_4=\sigma_3$. They imply $$(N+1)(\lambda_{1,0}-\lambda_{1,1})=\left(\sum_{p=0}^{N-2}\lambda_{1,p}-\lambda_{1,p+1}\right)+\lambda_4-\lambda_{1,0}+\lambda_{1,N-3}-\sigma_3=\lambda_4-\sigma_3=0
$$ so that $\lambda_{1,0}=\lambda_{1,1}$ and $\mu=0.$ Similarly, $\nu=0.$ Therefore the relations $(*_1)$ and $(*_2)$ now show that (PBW2') is satisfied, so that $A'$ is an algebra defined by a potential.

Thus if the characteristic of $\kk$ does not divide $(N+1)!,$ any PBW deformation of $A$ is defined by a potential.
This is not true however in arbitrary characteristic, as the following example shows. Let us fix $N=3$ and $\car(\kk)=2.$ Choose a PBW deformation $A'$ of $A$ defined by $\varphi_0=0$, $\varphi_1=0,$ and $\varphi_2(\partial_{a_1}W_4)=0,$ $\varphi_2(\partial_{a_2}W_4)=a_4a_3$, $\varphi_2(\partial_{a_3}W_4)=a_4a_1,$ $\varphi_2(\partial_{a_4}W_4)=a_1a_3.$ It is easy to check that this is indeed a PBW deformation of $A$ ((PBW1), (PBW3) and (PBW4) are obvious, and (PBW2) is equivalent to $(*_1)$ and $(*_2)$, which are simply the expressions of $\partial_{a_1}W_4$ and $\partial_{a_2}W_4$). However, (PBW2') does not hold, so that $A'$ is not defined by a potential (condition (PBW2') is a necessary condition for $A'$ to be defined by a potential, regardless of the characteristic of $\kk$).

Note that if in this last example we change $\varphi_2$ to: $\varphi_2(\partial_{a_1}W_4)=a_3a_4,$ $\varphi_2(\partial_{a_2}W_4)=0$, $\varphi_2(\partial_{a_3}W_4)=a_4a_1,$ $\varphi_2(\partial_{a_4}W_4)=a_1a_3,$ then we have a PBW deformation such that  (PBW2') is satisfied and that is defined by a potential ($W_3+W_2$ where $W_2=a_4a_1a_3$). Note that Theorem \ref{NSCpbw2'} does not apply here.

Finally, if in the previous example we assume that $\car(\kk)=3$ and  we again change $\varphi_2$ to: $\varphi_2(\partial_{a_1}W_4)=a_1^2,$ $\varphi_2(\partial_{a_2}W_4)=0$, $\varphi_2(\partial_{a_3}W_4)=0,$ $\varphi_2(\partial_{a_4}W_4)=0,$ then we have a PBW deformation such that  (PBW2') is satisfied, but it is not defined from a potential (the only candidate being given by $W_2=-a_1^3,$ but $\partial_{a_1}(W_2)=-3a_1^2=0\neq -\varphi_2(\partial_{a_1}W_4)$), thus showing that Theorem \ref{NSCpbw2'} is not applicable when the characteristic of $\kk$ divides $N!.$

\eex
%%%%%%%%%%%%%%%%%%%%%%%%%%%%%%%%%%%%%%%%%%%%%%%%%%%%%%%%%%%%%%%%%%%%%%%%%%%%%%%%%%%%%%%%%%%%%%%%%%%%%%%%%%%%%%%%%%%%%%%%%%%%%%%%%%%%%%%%%%%%%%%%%%%%%%%%%%%%%%%%%%%%%%%%%%%%%%%%%%%%%%%%%%%%%%%%%%%%%%%%%%%%%%%%%%%%%%%%%%%%%%%%%%%%%%%%%%%%%%%%%%%%%%%%%%%%%%%%%%%%%%%%%%%%%%%%%%%%%%%%%%%%%%%%%%%%%%%%%%%%%%%%%%%%%%%%%%%%%%%%%%%%%%%%%%%%%%%%%%%%%%%%%%%%%%%%%%%%%%%%%%%%%%%%%%%%%%%%%%%%%%%%%%%%%%%%
%%%%%%%%%%%%%%%%%%%%%%%%%%%%%%%%%%%%%%%%%%%%%%%%%%%%%%%%%%%%%%%%%%%%%%%%%%%%%%%%%%%%%%%%%%%%%%%%%%

%%%%%%%%%%%%%%%%%%%%%%%%%%%%%%%%%%%%%%%%%%%%%%%%%%%%%%%%%%%%%%%%%%%%%%%%%%%%%%%%%%%%%%%%%%%%%%%%%%
\flushleft{{\hrulefill\hspace*{10cm}}\\
\scriptsize{\sc  Roland Berger {\rm and} Rachel Taillefer \\Laboratoire de Math\'ematiques de l'Universit\'e de Saint-Etienne,\\Facult\'e des Sciences et Techniques,\\23 Rue Docteur Paul Michelon,\\42023 Saint-Etienne Cedex 2,\\France.\\ E-mail:} roland.berger@univ-st-etienne.fr and rachel.taillefer@univ-st-etienne.fr

\end{document}